The Hebrew University of Jerusalem

Faculty of Science

Einstein Institute of Mathematics

Master's Thesis in the subject of:

# Homological quantum error correcting codes and systolic freedom.

Submitted by: Ethan Fetaya

Guided by Prof. Alex Lubotzky

## Acknowledgements:

I would like to thank my advisor, Prof. Alex Lubotzky for help and guidance during my research and for teaching me how to write a scientific paper. I would also like to thank Prof. Józef Dodziuk for many fruitful discussions.

# Contents

# Introduction

This thesis is about homological quantum error correcting codes. These codes are best understood by looking first at classical error correcting codes.

**Definition 0.1.** *A classical linear code $C$, is a $k$ dimensional vector subspace of $\mathbb{Z}_2^n$. Its parameters are $[n, k, d]$ when $n$ and $k$ are the dimensions mentioned above, and the distance $d$ is the minimal Hamming distance between two elements of $C$. The Hamming distance is the number of coordinates on which two elements of $\mathbb{Z}_2^n$ differ (this defines a metric on $\mathbb{Z}_2^n$).*

We think of $C$ as the legitimate words to be transmitted over a noisy channel. The noise may changes the words we send, and we would want to reconstruct the original word after it has been altered. If $d = 2t + 1$ then the balls or radius $t$ in the Hamming distance around all the elements of $C$ are disjoint. This means that if the noise has changed at most $t$ coordinates then we can recover the message (by changing it to the closest element of $C$). Such a code is a $t$ error correcting code. The quality of a code is mainly measured by two parameters - rate and relative distance

**Definition 0.2.** *The rate $R$ of an $[n, k, d]$ code is $R = \frac{k}{n}$. The relative distance $\delta$ of a code is $\delta = \frac{d}{n}$. A good code is a family of codes $\{C_m\}_{m=1}^\infty$ such that $n \to \infty$ and the rates and relative distances are bounded away from zero.*

Since $n$ represents how many bits were physically transmitted, and $k$ represents how many bits of information were transmitted, the rate $R = \frac{k}{n}$ represents how efficient the code is in transmitting information. The relative distance $\delta = \frac{d}{n}$ tells us how fault tolerant the code is.

In quantum computation, the picture seems very different at first. The qubit, the quantum computer analog of the bit, is a vector in $\mathbb{C}^2$ of norm one, and quantum noise, unlike the classical "bit-flip", is some unitary transformation. Even though the task of quantum error correcting codes seems much more complicated, quantum error correcting codes can be constructed (and will be defined in Chapter 2) and are measured with the same parameters. A quantum code with parameters $[[n, k, d]]$ is a code where $n$ qubit are physically sent, $k$ qubits of information are sent, and if $d = 2t + 1$ then the code can correct $t$ errors (i.e. errors that effect at most t coordinates).

It turns out (for a detailed explanation see 2.2.3, which follows [24]) that the problem of constructing quantum error correcting codes, may lead to a combinatorial problem of a similar flavor as the classical codes. More specific, the CSS codes (called after A.M. Calderbank, P. Shor and A. Steane) are constructed from the following situation:

**Theorem 0.1.** *Let $V_1$ and $V_2$ be orthogonal subspaces of $\mathbb{Z}_2^n$ and let $W_j = V_j^\perp$ (so $V_j \subseteq W_{3-j}$). Then one can construct a quantum error correcting code (so called CSS codes) with parameters $[[n, k = dim(W_1/V_1) = dim(W_2/V_2), d]]$ when the distance $d$ is the minimal support of a vector orthogonal to one space, but is not in the other, i.e. the minimal Hamming weight of a vector in $(W_1/V_1) \cup (W_2/V_2)$.*

So construction of CSS codes boils down to finding for $n \to \infty$ orthogonal subspaces $V_1$ and $V_2$, and the code will have parameters $k$ and $d$ as described above.

One situation in which one gets orthogonal spaces in $\mathbb{Z}_2^n$ is the following: Let $M$ be a compact manifold of dimension $n$, and let $S$ be a triangulation of the manifold. For all $p \leq n$ one can define $Y_p$ as the $\mathbb{Z}_2$-space of all p-chains of $S$. i.e. the space spanned over $\mathbb{Z}_2$ by the p-simplices. This space has two well known subspaces, the boundaries $B_p = Im(\partial_{p+1})$ and the cycles $Z_p = Ker(\partial_p)$. These subspaces are not orthogonal, but since $B_p \subseteq Z_p$ then $B_p$ and $Z_p^\perp$ are orthogonal. We can define a CSS code by choosing $V_1 = B_p$ and $V_2 = Z_p^\perp$ both subspaces of $Y_p$, these codes will be called homological codes. To understand the distance of such codes we need to look at $Z_p/B_p$ and $B_p^\perp/Z_p^\perp$, the first is the $\mathbb{Z}_2$ homology is dimension p (with the basis defined by S), and the second can be taught of as the $\mathbb{Z}_2$ homology in dimension $n - p$ (with the basis defined by $S^*$ the dual cellulation).

**Definition 0.3.** *Let M be a connected compact manifold, with a triangulation $S$ and Assume $H_p(M, \mathbb{Z}_2) \neq 0$ for some natural number p. Let $S_p$ be the simplexes of dimension p, then it is a basis for p-chains with $\mathbb{Z}_2$ coefficients. The combinatorial systole of $M$ in dimension p, $csys_p(M, S)$, is defined as the minimal support of a $p-$cycle not homological to zero, in this basis given by $S_p$.*

We can conclude the properties of homological codes in the following theorem:

**Theorem 0.2.** *Let $M$ be a connected compact manifold, $S$ a triangulation and $S^*$ a dual triangulation (given by the poincaré duality). Assume $H_p(M, \mathbb{Z}_2) \neq 0$, then we can build a CSS code with $V_1 = B_p$ the boundaries, and $V_2 = Z_p^\perp$ the perpendicular space to the cycles. This is a quantum error correcting code with parameters $[[n, k, d]]$ when $n = |S_p|$, $k = dim(H_p(M, \mathbb{Z}_2))$, and $d = min\{csys_p(M, S), csys_{n-p}(M, S^*)\}$.*

The proof of these two theorems will be detailed in 2.2.4

In [26] Zémor showed some examples of such codes and conjectured

**Conjecture 0.1.** *(Zémor [26]) There exists some constant $C$, such that all homological quantum error correcting codes satisfy*

$$R\delta^2 \leq C \cdot n^{-2} \tag{1}$$

This conjecture implies that no good homological codes exist. The distance and the rate of these codes come from different geometric properties of the manifold. For this reason we look at a different inequality (that is equivalent to Zámors conjecture when looking at a constant manifold) that involves only the distance and not the rate

**Conjecture 0.2.** *There exists some constant $C$, such that all homological quantum error correcting codes satisfy*

$$d^2 \leq C \cdot n \tag{2}$$

Both conjectures are not true in the general case, and M.H. Freedman, D.A. Meyer and F. Luo in[21] construct a counterexample negating both. The work in [21] shows a strong correlation between conjecture 0.2 and systolic freedom.

**Definition 0.4.** *Let $M$ be a compact Riemannian manifold and assume $H_i(M, \mathbb{Z}_2) \neq 0$. Define the $\mathbb{Z}_2-$systole of $M$ as $sys_i(M, \mathbb{Z}_2) = \inf_{0 \neq \sigma \in H_i(M, \mathbb{Z}_2)} vol(\sigma)$, the minimal volume of a cycle not homological to zero in $\mathbb{Z}_2$ coefficients.*

**Definition 0.5.** *Let $n$, $p$ and $q$ be natural numbers such that $n = p+q$. The dimension $n$ has weak $(p,q)$ $\mathbb{Z}_2-$systolic freedom if*

$$\inf_M \frac{sys_n(M, \mathbb{Z}_2)}{sys_q(M, \mathbb{Z}_2) \cdot sys_p(M, \mathbb{Z}_2)} = 0 \tag{3}$$

*when the infimum is over all compact Riemannian manifolds of dimension $n$ with $H_p(M, \mathbb{Z}_2) \neq 0$.*

In [21] it was proven that in dimension 3 weak (2,1) $\mathbb{Z}_2-$systolic freedom exists, and that with systolic freedom one can construct codes that do not obey the $d^2 \leq C \cdot n$ bound. The proof of $\mathbb{Z}_2$-systolic freedom is by a beautiful construction, yet [21] has many missing details, and some small inaccuracies. In Chapter 3 we present a more detailed and accessible proof of the construction in [21].

In Chapter 4, which is the main novelty of this thesis, we use the fact that in dimension 2 there is no systolic freedom, i.e. it has systolic rigidity, to prove that conjecture 0.2 holds for all homological codes from surfaces.

**Theorem 0.3.** *There exists a constant $C$ such that every homological code from a surface has distance bounded by $d^2 \leq C \cdot n$*

The work here and in [21] show a strong correlation between systolic freedom and the $d^2 \leq C \cdot n$ bound. It is not yet known if good homological quantum error correction codes exist, but the search for these codes is connected to the search of dimensions with "stronger" freedom, and one can only hope to do this in dimension larger then 2.

# Notation and convention

We will use the Bachmann-Landau asymptotic notation. We will say that

i. $f = \mathcal{O}(g)$, if there exists constants $C$ and $x_0$ such that $f(x) \leq C \cdot g(x)$ for all $x > x_0$.

ii. $f = \Omega(g)$, if there exists constants $C$ and $x_0$ such that $f(x) \geq C \cdot g(x)$ for all $x > x_0$.

iii. $f = \Theta(g)$, if $f = \mathcal{O}(g)$ and $f = \Omega(g)$.

iv. $f = o(g)$ if $lim_{x \to \infty} \frac{f(x)}{g(x)} = 0$

I will also use in the first two chapters the bra-ket notation normally used in QM (quantum mechanics). Let $\mathcal{H}$ be a complex inner product space. A vector $\psi$ in $\mathcal{H}$ will be called a *ket* and written $|\psi\rangle$. The dual vector to $|\psi\rangle$ will be called *bra* and written $\langle\psi|$. The inner product $(|\psi\rangle,|\varphi\rangle)$ will be written $\langle\psi|\varphi\rangle$ and can be thought of as the functional $\langle\psi|$ operating on the vector $|\varphi\rangle$: $\langle\psi|(|\varphi\rangle)$. A notable deviation from standard mathematical convention is that in this convention, the inner product is linear with regard to the right argument and not the left argument.

If A is a linear operator, then the dual vector to $A|\psi\rangle$ is $\langle\psi|A^\dagger$, when $A^\dagger$ is the adjoint operator to $A$. We will understand $\langle\varphi|A|\psi\rangle$ as $(|\varphi\rangle, A|\psi\rangle)$ or dually $(A^\dagger|\varphi\rangle, |\psi\rangle)$. Let $\lambda$ be an eigenvalue of A, then we will sometimes use the same notation for the specific eigenvector we are working with i.e. $A|\lambda\rangle = \lambda|\lambda\rangle$. When there is degeneracy and the dimension of the eigenspace is greater than one, we will add another index. We will define the outer product of a ket $|\psi\rangle$ and a bra $\langle\varphi|$ as an operator marked $|\psi\rangle\langle\varphi|$. The action of $|\psi\rangle\langle\varphi|$ on a vector $|\alpha\rangle$ is defined as $|\psi\rangle\langle\varphi|(|\alpha\rangle) = |\psi\rangle\langle\varphi|\alpha\rangle = \langle\varphi|\alpha\rangle|\psi\rangle$.

The tensor product of two (or more) vectors $|\psi\rangle \otimes |\varphi\rangle$ will be written in short form as $|\psi\rangle|\varphi\rangle$. In quantum computation we will work in a multiple tensor product of 2 dimensional spaces. We will mark the base of the 2 dimensional space by $|0\rangle$ and $|1\rangle$, and $|e_1...e_n\rangle$ will mean $|e_1\rangle \otimes ... \otimes |e_n\rangle$ when $e_i \in \{0, 1\}$.

# Chapter 1

# Introduction to Quantum Mechanics

This chapter gives a short review of the foundations of quantum mechanics. Anyone familiar with the subject can freely skip this chapter.

Quantum mechanics are the laws that determine the behavior of small scale objects. The postulates of quantum mechanics are fundamentally different from those of classical mechanics, and therefore a quantum computer in which the basic unit (the qubit) is governed by the laws of quantum mechanics has different abilities and limitations than those of a regular computer.

## 1.1 The Postulates of quantum mechanics

Unlike classical mechanics which has $\mathbb{R}^n$ as its natural setting, quantum mechanics' natural settings are complex Hilbert spaces.

**Postulate 1.** *To each closed quantum mechanical system there is a Hilbert space $\mathcal{H}$ called the* state space. *The unit vectors in this state space completely define the system.*

**A few notes:** First, a closed system is a system without any interaction with its surroundings. Such systems (unless you take the whole universe as your system) do not really occur, but normally the systems physicists work with, are approximately closed. The influence of the environment is precisely what causes the *noise* whose influence quantum error correcting codes try to reverse. Second, two vectors in $\mathcal{H}$ that differ by a scalar factor represent the same physical system so we limit ourselves to unit vectors without loss of generality.

In most physical systems, the state space will be infinite dimensional, but in quantum computation only finite dimensional systems are considered. The most basic unit of a quantum computer is the qubit. The qubit is a physical system whose state space is two dimensional like the electron spin for example. The vectors in an orthonormal basis for this space are usually marked $|0\rangle$ and $|1\rangle$.

The main difference between a normal computer bit and the qubit is that a bit can only be in one of two states - zero or one, while the qubit can be in *superposition* of these two states i.e. can be in any linear combination of this two states. The general qubit can be written as $|\psi\rangle = a|0\rangle + b|1\rangle$ where $a, b \in \mathbb{C}$ satisfy $|a|^2 + |b|^2 = 1$.

**Time Evolution**

**Postulate 2.** *The evolution of a closed quantum system in time is given by a* unitary *operator. Moreover, to each quantum system there is a Hermitian operator H called the Hamiltonian, and the evolution of the system is given by the Schrödinger equation $i\hbar \frac{d|\psi\rangle}{dt} = H|\psi\rangle$ where $\hbar$ is a physical constant.*

When the Hamiltonian is constant in time it is easy to see that $|\psi(t)\rangle = exp\left[\frac{-i(t-t_0)H}{\hbar}\right]|\psi(t_0)\rangle$ is the solution to the Schrödinger equation with initial state $|\psi(t_0)\rangle$ at $t = t_0$. From this it is clear that the time evolution is given by the operator $U(t, t_0) = exp\left[\frac{-i(t-t_0)H}{\hbar}\right]$ and that $U(t, t_0)$ is indeed unitary since H is hermitian.

The fact that H is an Hermitian operator implies that it has a complete orthonormal system of eigenvalues, therefore $H = \sum_E E|E\rangle\langle E|$ where $H|E\rangle = E|E\rangle$. The eigenstates (i.e. eigenvectors) of H are called stationary states. This name comes from the fact that if we take $|\psi(t_0)\rangle = |E\rangle$ then $|\psi(t)\rangle = exp\left[\frac{-i(t-t_0)H}{\hbar}\right]|E\rangle = exp\left[\frac{-i(t-t_0)E}{\hbar}\right]|E\rangle$ this means that $|\psi(t)\rangle$ and $|\psi(t_0)\rangle = |E\rangle$ differ by a scalar factor and therefore represent the same physical system. A linear combination of the stationary states, on the other hand, will not in general stay stationary since the phase factor (the scalar) changes differently for each of its components.

### Measurements

The notion of measurements is the most significant conceptual difference between quantum and classical mechanics. The measurement postulate is defined a bit differently between quantum computation books [23] and physicists' quantum mechanics books, with the quantum computation books giving a slightly broader definition. We will cover both approaches here, beginning with the usual QM approach.

**Postulate 3.** *Every measurable physical quantity is described by an operator A operating on the space $\mathcal{H}$. The possible results of the measurement are the eigenvalues of A. If $\lambda$ is an eigenvalue of A then the probability of measuring $\lambda$ when the system is in the state $|\psi\rangle$ is $\langle\psi| P_\lambda^\dagger P_\lambda |\psi\rangle$, when $P_\lambda$ is the projection operator onto the $\lambda$ eigenspace. After the measurement, the state changes into $\frac{1}{\sqrt{\langle\psi|P_\lambda^\dagger P_\lambda|\psi\rangle}} P_\lambda |\psi\rangle$ when $\frac{1}{\sqrt{\langle\psi|P_\lambda^\dagger P_\lambda|\psi\rangle}}$ is the normalization factor.*

Since the result of a physical measurement is a real number, all the eigenvalues of A must be real. If we assume that A has a complete orthonormal system, we will conclude that A is Hermitian. For example, let the system be in a state $|\psi\rangle$ and measure the Hamiltonian H (the physical interpretation of this measurement is the energy of the system). If there is no degeneracy, i.e. all eigenspaces are of dimension 1, then $|\psi\rangle = \sum_{E_n} a_n |E_n\rangle$ where $\sum_n |a_n|^2 = 1$, then when we measure H we get $|E_n\rangle$ with probability $|a_n|^2$. The system after we measured $E_n$ is in state $|E_n\rangle$

The measurement postulate has two serious impacts on the behavior of quantum systems. First, it has probability as an intrinsic aspect. If you have a state $|\psi\rangle = \frac{1}{\sqrt{2}} |0\rangle + \frac{1}{\sqrt{2}} |1\rangle$ and you measure whether it is in the $|0\rangle$ state or the $|1\rangle$ state, the theory states that the result is determined randomly with equal probabilities. This means that two similar states can give different results during the same measurement. This contradicts the classical notion of *determinism*. The fact that QM is inherently probabilistic bothered many physicists, notably Albert Einstein who famously said "God does not play dice with the universe". His attempts to disprove this probabilistic interpretation brought him to design the "EPR paradox"[1] which leads to *quantum entanglement*. Entanglement will play a crucial role in quantum computation. The second important repercussion of this postulate is the fact that you change the state of the system by measuring it. This adds another level of complexity in designing quantum algorithms - every measurement can destroy data in an irreversible way. One of the known properties of quantum measurements is the uncertainty principle which we will cover later on.

The measurement postulate given in [23] defines measurement as a collection $\{M_\lambda\}$ of *measurable operators* satisfying $\sum_\lambda M_\lambda^\dagger M_\lambda = I$. The index $\lambda$ is the outcome of the measurement, and the probability of measuring $\lambda$ is $\langle\psi|\,M_\lambda^\dagger M_\lambda\,|\psi\rangle$ therefore the demand $\sum_\lambda M_\lambda^\dagger M_\lambda = I$ guarantees that the probabilities sum up to one. The state $|\psi\rangle$ after the measurement is $\frac{M_\lambda|\psi\rangle}{\sqrt{\langle\psi|M_\lambda^\dagger M_\lambda|\psi\rangle}}$. If we take $M_\lambda$ to be the projection operator into the $\lambda$ eigenspace of A, we will get the previous definition. This kind of measurement will be called *projective measurement*.

**Composite systems**

We know how to describe a physical system, for example the qubit. The next postulate gives us the description of a system that is the composition of several smaller systems

**Postulate 4.** *If the Hilbert spaces $\mathcal{H}_1, ..., \mathcal{H}_n$ represent separate components of a composite physical system, then the Hilbert space representing the composite system $\mathcal{H}$ is the tensor product of component spaces i.e. $\mathcal{H} = \mathcal{H}_1 \otimes ... \otimes \mathcal{H}_n$.*

Since the basic unit we are working with, the qubit, is represented by $\mathbb{C}^2$, the n-qubit system is represented by $\mathbb{C}^{2^n}$. If $|0\rangle$ and $|1\rangle$ are the basis of each qubit space, then $|e_1...e_n\rangle\, e_i \in 0,1$ is the basis of the composite system.

An interesting consequence of this postulate is *Quantum entanglement*. Let us think of the 2-qubit system $\mathbb{C}^2 \otimes \mathbb{C}^2$. If a state $|\psi\rangle$ of the composite system is indecomposable, then we cannot think of one qubit being in one state and one being in another - they are entangled. The most well known example of this is the "EPR pair" $|\psi\rangle = \frac{1}{\sqrt{2}}(|01\rangle + |10\rangle)$ (this state is easily realized by a physical system [1]). In this system, when you measure the state of one qubit, you instantly change the state of the second qubit into the opposite state. If we have two electrons that are an EPR pair, then when you measure one of the electrons, you change the other and that change must be instantaneous. This instantaneous change contradicts the idea of *locality* (that information cannot travel faster than the speed of light) and this is the EPR "paradox". This entanglement, however, was proven experimentally, and it can be proven by Bell's inequality that no physical theory that explains quantum phenomenons can be local and deterministic [2].

## 1.2 Measurement expectation and the uncertainty principle

The idea behind Heisenberg's uncertainty principle is that if two measurement operators do not commute, then you cannot know both well. The most common example is the location operator (on a certain axis) and momentum operator (on that same axis). The more you know one, the less you know about the other. Before we prove the uncertainty principle, we must explain what "knowing something well" means.

Since the measurement outcome is determined randomly, the smaller the standard deviation, the closer you probably are to the expected value, and the outcome is better known. In the case of projective measurements, there is a simple formula for the expected value and the variance. It is given in the following lemma

**lemma 1.1.** *If M is a projective measurement, then when measuring M on the state $|\psi\rangle$ the expected value is given by $E(M) = \langle\psi| M |\psi\rangle$ and the variance by $\sigma^2(M) = E((M - E(M))^2) = \langle\psi| (M - E(M))^2 |\psi\rangle$.*

**Proof.** *$E(M) = \Sigma_m mp(m) = \Sigma_m m \langle\psi| P_m |\psi\rangle = \langle\psi| \Sigma_m m P_m |\psi\rangle = \langle\psi| M |\psi\rangle$. The variance follows immediately.*

Now we are ready to state and prove the uncertainty principle

**The Heisenberg uncertainty principle.** *Let A and B be Hermitian operators, and $|\psi\rangle$ be a quantum state. Then $\sigma(A)\sigma(B) \geq \frac{|\langle\psi|[A,B]|\psi\rangle|}{2}$.*

**Proof.** *For any two Hermitian operators C,D and state $|\psi\rangle$, $\langle\psi| CD |\psi\rangle = x + iy$ for some $x, y \in \mathbb{R}$. It is easy to see that $\langle\psi| [C, D] |\psi\rangle = 2iy$ and $\langle\psi| \{C, D\} |\psi\rangle = 2x$ when $[A, B] = AB - BA$ is the commutator and $\{A, B\} = AB + BA$ is the anti-commutator.
From this we get*

$$|\langle\psi| [C, D] |\psi\rangle|^2 + |\langle\psi| \{C, D\} |\psi\rangle|^2 = 4|\langle\psi| CD |\psi\rangle|^2. \qquad (1.1)$$

*Using the Cauchy-Schwarz inequality we get*

$$|\langle\psi|CD|\psi\rangle|^2 = |\langle C^\dagger\psi|D\psi\rangle|^2 \leq \langle\psi|C^2|\psi\rangle\langle\psi|D^2|\psi\rangle. \qquad (1.2)$$

*If we combine the (1.1) with (1.2), and omit the positive anti-commutator part to get*

$$\frac{1}{4}|\langle\psi|[C,D]|\psi\rangle|^2 \leq \langle\psi|C^2|\psi\rangle\langle\psi|D^2|\psi\rangle. \qquad (1.3)$$

*If we substitute C=A-E(A), D=B-E(B) and take the square root, we get the uncertainty principle.* □

If we look at the other case, where A and B commute, then we can mutually diagonalize them and get an orthonormal base $|n,m,k\rangle$ where $A|n,m,k\rangle = n|n,m,k\rangle$, $B|n,m,k\rangle = m|n,m,k\rangle$ (the extra index k is needed when there is degeneracy in the eigenspace). The state $|n,m,k\rangle$ has zero standard deviation when measuring A or B (since the outcome is n or m with probability 1) and therefore when A and B commute you can measure both together as accurately as you want.

A fascinating way to see the physical ramifications of the uncertainty principle is the *Stern-Gerlach experiment* explained throughly in [23].

## 1.3   The density operator

There is an alternative way to formulate quantum mechanics through the *density operator*. It is not as intuitive as the previous formulation but it has two advantages. First, it lets us work with ensembles of quantum states. Second, it lets us describe subsystems of a composite quantum system. The second property is very important in understanding quantum noise, as the noise comes from the fact that our system is not closed and we need to think of it as part of a larger system within its environment.

**Definition 1.1.** *Let $\{|\psi_i\rangle\}_{i=1}^n$ be quantum states, If the system is in the state $|\psi_i\rangle$ with probability $p_i$, where $\sum_i p_i = 1$ then we call $\{p_i, |\psi_i\rangle\}$ an ensemble of pure states. The density operator for this ensemble is defined as $\rho = \sum_i p_i |\psi_i\rangle\langle\psi_i|$.*

It is very important to note the difference between a pure state (i.e. a vector in the state space) and an ensemble. The pure state $|\psi\rangle = \frac{1}{\sqrt{2}}(|0\rangle+|1\rangle)$ is different from an ensemble that is $|0\rangle$ with probability $\frac{1}{2}$ and $|1\rangle$ with probability $\frac{1}{2}$, even though the probability to measure if the quantum state is in $|0\rangle$ or $|1\rangle$ is identical. We will show how the postulates of QM can be stated with the density operator formalism.

**Theorem.** *Let $\rho$ be the density operator of a system. If the time evolution of the system is given by $U$, then the evolution of the density operator is given by $U\rho U^\dagger$. If a measurement is performed by the operators $\{M_m\}$, then the probability of measuring $m$ is $p(m) = tr(M_m^\dagger M_m \rho)$ and the density operator after the measurement is $\frac{M_m \rho M_m^\dagger}{tr(M_m^\dagger M_m \rho)}$ when $m$ is the outcome of the measurement.*

**Proof.** *If the evolution of the system is given by $U$, then each state $|\psi_i\rangle$ transforms to $U|\psi_i\rangle$. For the density operator we get*

$$\rho = \sum_i p_i |\psi_i\rangle \langle\psi_i| \to \sum_i p_i U |\psi_i\rangle \langle\psi_i| U^\dagger = U\rho U^\dagger \tag{1.4}$$

*To prove the statement we will use a useful and easy to check trace formula*

$$tr(A |\psi\rangle \langle\psi|) = \langle\psi| A |\psi\rangle \tag{1.5}$$

*now to calculate $p(m)$ from the law of total probability*

$$p(m) = \sum_i p_i p(m|i) = \sum_i p_i \langle\psi_i| M_m^\dagger M_m |\psi_i\rangle = \sum_i p_i tr(M_m^\dagger M_m |\psi_i\rangle \langle\psi_i|) = tr(M_m^\dagger M_m \rho) \tag{1.6}$$

*If the result of the measure is $m$, then the state $|\psi_i\rangle$ after the measurement is $|\psi_i^m\rangle = \frac{M_m |\psi_i\rangle}{\sqrt{\langle\psi_i| M_m^\dagger M_m |\psi_i\rangle}} = \frac{1}{\sqrt{p(m|i)}} M_m |\psi_i\rangle$. The density operator after the measurement will be*

$$\rho_m = \sum_i p(i|m) |\psi_i^m\rangle \langle\psi_i^m| = \sum_i \frac{p(i|m)}{p(m|i)} M_m |\psi_i\rangle \langle\psi_i| M_m^\dagger \tag{1.7}$$

*From Bayes formula we drive the desired equality*

$$\rho_m = \sum_i \frac{p_i}{p(m)} M_m |\psi_i\rangle \langle\psi_i| M_m^\dagger = \frac{M_m \rho M_m^\dagger}{tr(M_m^\dagger M_m \rho)} \tag{1.8}$$

□

We will cite a characterization of the density operator from [23, p. 101] without proof, as it is easy to verify.

**Theorem.** *An operator $\rho$ is the density operator for some ensemble iff, $\rho$ is a positive operator with trace equal to one. A density operator $\rho$ is a density operator of a pure state (i.e. $\rho = |\psi\rangle \langle\psi|$ for some state $|\psi\rangle$) iff $tr(\rho^2) = 1$*

**The reduced density operator**

If the system is composed of two subsystems A and B, described by a density operator $\rho^{AB}$, then the density operator that describes the A component is the partial trace over B, $\rho^A = tr_B(\rho^{AB})$. It can be shown [23, p. 107] that this definition yields the desired time evolution and measurement. In the simple case of $\rho^{AB} = \rho \otimes \sigma$ this definition does give $\rho^A = \rho$ as expected.

The connection to quantum noise is immediate - if one starts from a open system A, then it can be thought of as a part of the composite system containing it and its environment E. The time evolution of the composite system is given by a unitary transformation, so $\rho^{AE} \to U\rho^{AE}U^\dagger$. The original system transforms to $tr_E(U\rho^{AE}U^\dagger)$, this transformation does not have to be unitary and in fact in the next example we will show that even if U does not change A, $\rho^A$ might change.

**Example:** Let A be in the pure state $\frac{1}{\sqrt{2}}|0\rangle + \frac{1}{\sqrt{2}}|1\rangle$ and B be in the state $|0\rangle$. The composite system is in the state $\frac{1}{\sqrt{2}}|00\rangle + \frac{1}{\sqrt{2}}|10\rangle$.
Let U be the controlled-not transformation

$$U|00\rangle = |00\rangle \quad U|01\rangle = |01\rangle \quad U|10\rangle = |11\rangle \quad U|11\rangle = |10\rangle \qquad (1.9)$$

The operator U, of course, does not change A, but the composite state transforms to $\frac{1}{\sqrt{2}}|00\rangle + \frac{1}{\sqrt{2}}|11\rangle$ and a simple calculation shows the reduced density matrix after the transformation is $\rho_A = \frac{1}{2}|0\rangle\langle 0| + \frac{1}{2}|1\rangle\langle 1|$. The system is now not only changed, it is not in a pure state anymore.

# Chapter 2

# Quantum Computation and Information

With the QM language and foundations given, we can now describe the *quantum circuit model* for a quantum computer, which will be the model we will use for quantum computation. We will explain later the basics of quantum error correction, and formulate the homological error correcting codes that is the heart of this thesis.

## 2.1 Quantum circuit model

The quantum computer in the quantum circuit model is made out of four parts: qubits, gates, measurements and a classical part. Qubits and measurement have been explained previously, and the classic computer is obvious (it is also an avoidable part of this model) which leaves gates to be described.

**Quantum gates:** The quantum gate is an operator that works on an n-qubit state space. Since time evolution of a system is unitary, a gate is an unitary operator. We will assume that every operator can be implemented perfectly (for a treatment of approximating quantum gates see [23]). Every gate is reversible, unlike the classical computer where most of the basic gates such as *and* and *or* are irreversible. This problem can be corrected and every classical operation can be emulated on a quantum computer. For example the *xor* gate cannot be implemented as is, since it is two to one and as such is not reversible, but it can be emulated be the gate $U\left|x_1, x_2\right> = \left|x_1, x_1 + x_2 (mod 2)\right>$ on the computational basis described in 1.1.1.

The Pauli matrices X,Y,Z and the Hadamard gate H are important gates, all operate on a single qubit, defined as following:

$$X = \begin{pmatrix} 0 & 1 \\ 1 & 0 \end{pmatrix} \quad Y = \begin{pmatrix} 0 & -i \\ i & 0 \end{pmatrix} \quad Z = \begin{pmatrix} 1 & 0 \\ 0 & -1 \end{pmatrix} \qquad (2.1)$$

$$H = \frac{1}{\sqrt{2}} \begin{pmatrix} 1 & 1 \\ 1 & -1 \end{pmatrix} \qquad (2.2)$$

The X gate is the quantum equivalent of a "bit-flip" transformation, while Z which is called the "phase-flip" has no classical analog. The operator Y can be thought of as the product of a bit-flip and phase-flip, since $ZX = iY$ and the $i$ is a global phase factor that does nothing to the state. The Hadamard gate is the discrete Fourier transform.

A summery of some of the important properties of the Pauli matrices:

**Theorem 2.1.** *The Pauli matrices $X, Y, Z$ satisfy:*

i. $\{X, Y, Z, I\}$ form a basis for $M_2(\mathbb{C})$

ii. $X^2 = Y^2 = Z^2 = I$.

iii. $det(X) = det(Y) = det(Z) = -1 \quad Tr(X) = Tr(Y) = Tr(Z) = 0$.

iv. $XY = iZ \quad YZ = iX \quad ZX = iY$.

v. $[X, Y] = 2iZ \quad [Y, Z] = 2iZ \quad [Z, X] = 2iY$.

**Quantum algorithms**

A quantum algorithm in this model is a circuit built out of the previously discussed parts. It is best understood by examples and we will show a few that illustrate some of the power and some of the problems with the quantum computer.

**Quantum parallelism:** Lets us assume that there is a function $f : \{0,1\}^n \to \{0,1\}$ that you want to compute by a quantum computer. Since the transformation $U|x\rangle = |f(x)\rangle$ ($x \in \{0,1\}^n$ describes a general basis element not an arbitrary state) is not (in general) reversible, it cannot be implemented as a gate. The operator $U_f|x,y\rangle = |x, y \oplus f(x)\rangle$ (where $\oplus$ means addition modulo 2) when $U_f : \mathbb{C}^{2^{n+1}} \to \mathbb{C}^{2^{n+1}}$, on the other hand, is unitary and can be implemented as a gate. To get quantum parallelism, we will use an important property of the Hadamard gate: $H^{\otimes n}|0...0\rangle = \left(\frac{|0\rangle+|1\rangle}{\sqrt{2}}\right) \otimes ... \otimes \left(\frac{|0\rangle+|1\rangle}{\sqrt{2}}\right) = \frac{1}{\sqrt{2^n}} \sum_{x \in \{0,1\}^n} |x\rangle$ i.e. the image of the state $|0...0\rangle$ is the sum of all basis elements with equal coefficient. This is a special case of a more general formula

$$H^{\otimes n}|x\rangle = \frac{1}{\sqrt{2^n}} \sum_{z \in \{0,1\}^n} (-1)^{z \cdot x} |z\rangle. \tag{2.3}$$

When $z \cdot x$ is the scalar product $\sum_i z_i x_i$.

Now we start with the n+1 qubit $|\psi_1\rangle = |0...0\rangle \otimes |0\rangle$ and pass it through a $H^{\otimes n} \otimes I$ gate to get $|\psi_2\rangle = H^{\otimes n} \otimes I |\psi_1\rangle = \frac{1}{\sqrt{2^n}} \sum_x |x\rangle \otimes |0\rangle$. In the next step we will use the $U_f$ gate on $|\psi_2\rangle$ to receive $|\psi_3\rangle = U_f |\psi_2\rangle = \frac{1}{\sqrt{2^n}} \sum_x |x\rangle \otimes |0 \oplus f(x)\rangle = \frac{1}{\sqrt{2^n}} \sum_x |x\rangle \otimes |f(x)\rangle$. In fact we computed all $2^n$ values of $f(x)$ by using the $U_f$ gate only once! The problem is that in order to receive this information, you have to measure the state and destroy most of the data. If you measure the value of $x$, you will get some answer $x_i$ randomly and the state after the measurement will be $|x_i\rangle \otimes |f(x_i)\rangle$. De facto, not only do we get only one answer, we get the answer for some random $x$ value. The power of quantum computation comes from the ability of some algorithms to use this global information without using measurements. The two best known examples of quantum algorithms that surpass their classical version, are the quantum Fourier transform [23] and Shor's prime factoring algorithm [25]. We will give a simple example, that might shed some light on how this could be done, with the *Deutsch-Jozsa algorithm*.

**The Deutsch-Jozsa algorithm:** This algorithm is a solution to the following game. Bob has a function $f : \{0,1\}^n \to \{0,1\}$ that is either the constant function, or a function that returns 1 and 0 an equal amount of times (then $f$ is said to be balanced). Alice can give Bob any $x \in \{0,1\}^n$ and he will compute $f(x)$ and return it to her. The question is, how many times must Alice and Bob exchange information in order for her to be absolutely certain what kind of function Bob has. The classical worst case scenario is of course $2^{n-1} + 1$. We will show that Alice can do it with only a single qubit exchanged between her and Bob.

Alice will start with the state $|\psi_1\rangle = |0\rangle^{\otimes n} |1\rangle$ and pass it through a $H^{\otimes(n+1)}$ gate to get $|\psi_2\rangle = \frac{1}{\sqrt{2^n}} \sum_x |x\rangle \left(\frac{|0\rangle - |1\rangle}{\sqrt{2}}\right)$. Next, Alice will send $|\psi_2\rangle$ to Bob who will execute the $U_f$ gate on $|\psi_2\rangle$. The first n qubits stay unchanged, but the last qubit after $U_f$ is performed is $|0 \oplus f(x)\rangle - |1 \oplus f(x)\rangle = (-1)^{f(x)}(|0\rangle - |1\rangle)$. So $|\psi_3\rangle = U_f |\psi_2\rangle = \frac{1}{\sqrt{2^n}} \sum_x (-1)^{f(x)} |x\rangle \left(\frac{|0\rangle - |1\rangle}{\sqrt{2}}\right)$. Alice receives $|\psi_3\rangle$ from Bob and uses a $H^{\otimes n} \otimes I$ gate on it. From eq. 2.3 we get

$$|\psi_4\rangle = \sum_z \sum_x \frac{1}{2^n} (-1)^{f(x)+x \cdot z} |z\rangle \left(\frac{|0\rangle - |1\rangle}{\sqrt{2}}\right) \quad (2.4)$$

For the last part of the algorithm, Alice measures the first n bits. The meaning of the measurement is that we take the operator $A|x\rangle = x|x\rangle$ when we think of x as a number given in base 2, to be the measured operator. The measurement postulate stated that the result is just the x (when we think of measuring x or $|x\rangle$ as the same thing) with probability that is equal to the square of that basis element's coefficient. In this case the probability of measuring $|0\rangle^{\otimes n}$ is just square of $\sum_x \frac{(-1)^{f(x)}}{2^n}$. When $f$ is constant,a we will get $\pm 1$ and the probability of measuring $|0\rangle^{\otimes n}$ is 1, when $f$ is balanced the $+1$ and $-1$ will cancel each other and the probability will be 0. So Alice can now know with probability one what type of function Bob is using after passing just a *single* qubit between her and Bob (and after using only 2 gates so there are no computational complexity problems).

Perhaps this algorithm is not the most useful, but it shows the flavor and abilities of a quantum computer without being too complicated.

## 2.2 Quantum Error Correction

When facing the problem of quantum error correction, there are a few fundamental problems that separate it from the classical counterpart.

- The error is described as a unitary operator and as such, there is a continuem of possible errors, unlike the finite classical case.

- Measurements can destroy data.

- Quantum gates are only approximated (since there is a continuem of possible gates) and therefore there are errors that come from this imperfect implementation.

We will ignore the last problem and assume that computations can be done perfectly (the problem of fault-tolerant computation is dealt in depth in [23]). Even when considering only the first two problems, the classical and quantum errors seem very different, so it seems surprising when CSS codes (that will be described later on) build a QECC (quantum error correction code) from two CECC (classical error correction code).

### 2.2.1 Quantum noise

Quantum noise comes from a unitary transformation that operates on the system and its environment as described in 1.3. If the system is in state $|\psi\rangle$, we can assume the composite system is in the state $|\psi\rangle \otimes |0\rangle_{env}$ when $|0\rangle_{env}$ describes the state of the environment in some, not necessary finite dimensional, space (this can be assumed without loss of generality by purification [23]). This composite system undergoes a unitary transformation $U$ that we want to reverse. In order to describe $U$ we will use the Pauli group.

**Definition 2.1.** *The* Pauli group *is the group $G_n$ of $2^{2n+1}$ operators $\{\pm X, \pm Z, \pm iY, \pm I\}^{\otimes n}$. The center of this group is $\{\pm I^{\otimes n}\}$, and we will mark $\bar{G}_n = G_n/\{\pm I^{\otimes n}\} = \{X, Z, iY, I\}^{\otimes n}$. Each element of these groups is an operator that is a tensor product (signed or unsigned) of Pauli matrices and the identity.*

We will use the fact that this is a group later on, at this moment we will look at $\bar{G}_n$ as a basis for $M_{2^n}(\mathbb{C})$, the space of operators on a n qubit space. In some of the literature (e.g. [24]), $\bar{Y} \equiv iY = ZX = \begin{pmatrix} 0 & 1 \\ -1 & 0 \end{pmatrix}$ is defined as the $Y$ Pauli matrix. This may be more convenient when discussing QEC but we lose the nice symmetrical properties of the Pauli matrices.

As previously stated what we are aiming at is to describe the action $U|\psi\rangle \otimes |0\rangle_{env}$. This is given by the following theorem

**Theorem.** *Let $\{E_a\}$ be the elements of $\bar{G}_n$. For any unitary operator on the composite system there are states $\{|e_a\rangle_{env}\}$* **generally not normalized or orthogonal** *such that $U |\psi\rangle \otimes |0\rangle_{env} = \sum_a E_a |\psi\rangle \otimes |e_a\rangle_{env}$.*

This can be easily proven directly for $n = 1$ and for larger $n$ by induction (see [24, Ch. 7]).

The idea behind QECC is that the elements of $\bar{G}_n$ represent all possible errors. Like in the case of CECC we cannot correct all possible errors only errors with weight beneath a certain threshold, and we assume that the larger errors are with negligible probability. The weight of $E_a \in \bar{G}_n$ is defined as the number of qubits on which $E_a$ acted non-trivially. We now take a subset $\mathcal{E} \subseteq \bar{G}_n$ and these are the errors we are trying to correct. The idea behind QEC is that if we start with a state $|\psi\rangle$, if we can measure what $E_a \in \mathcal{E}$ occurred then after the measurement you get $E_a |\psi\rangle \otimes |e_a\rangle_{env}$ and the affect of $E_a$ can be reversed since we know $E_a$ and it is unitary.
Usually we will take $\mathcal{E}$ to be the set of all the errors with weight$\leq t$ for some t. These are going to be the errors in $\bar{G}_n$ that we will want to fix.

### 2.2.2 Quantum Error Correcting codes

A quantum t-error correcting code will be a k-qubit state space embedded in an n-qubit state space such that it can correct $\mathcal{E}_t$ - all the errors with weight up to t. The next theorem states sufficient and necessary conditions for a code to correct t errors.

**Theorem.** *Let $\{|\bar{i}\rangle\}$ be a basis for the code space (the subspace that we are encoding). The code can correct errors in $\mathcal{E}$ iff for any $E_a, E_b \in \mathcal{E}$, $\langle \bar{j}| E_a^\dagger E_b |\bar{i}\rangle = C_{ab} \delta_{ij}$ when $C_{ab} = \langle \bar{i}| E_a^\dagger E_b |\bar{i}\rangle$ is* **independent of** $|\bar{i}\rangle$.

It is obvious that $\langle \bar{j}| E_a^\dagger E_b |\bar{i}\rangle = 0$ for $i \neq j$ is a necessary condition, because otherwise the states $E_b |\bar{i}\rangle$ and $E_a |\bar{j}\rangle$ cannot be distinguished by a measurement - we cannot know if we started with $|\bar{i}\rangle$ and the noise was $E_b$, or we started with $|\bar{j}\rangle$ and the noise was $E_a$. It is also clear that $\langle \bar{j}| E_a^\dagger E_b |\bar{i}\rangle = \delta_{ab} \delta_{ij}$ is sufficient, since if $\mathcal{H}$ is the code space then the spaces $\mathcal{H}_a = E_a \mathcal{H}$ are orthogonal and thus we can measure in what $\mathcal{H}_a$ we are in. This kind of code is called *nondegenerate*. The intuition behind the more general case (i.e. $\langle \bar{i}| E_a^\dagger E_b |\bar{i}\rangle$ is independent of $|\bar{i}\rangle$ but not $\delta_{ab}$) , is that if $\langle \bar{i}| E_a^\dagger E_b |\bar{i}\rangle$ does depends on $|\bar{i}\rangle$ we get some knowledge on the state when we measure $\mathcal{H}_a$ and

by so doing we inevitably change it, so in order to correct the code $C_{ab}$ must not depend on $|\bar{i}\rangle$ (for a full proof see [24, ch. 7]).

**Definition 2.2.** *The* distance *of a code C will be the minimal d such that there is an operator $E_a \in \bar{G}_n$ of weight d satisfying $\langle \bar{i}| E_a |\bar{j}\rangle \neq C_a \delta_{ij}$. Like in CECC, if d=2t+1 then the code can correct t errors (since $\bar{G}_n$ is a group). A k-qubit code in an n-qubit state space that has distance d will be marked [[n,k,d]]. It is a $2^k$ dimensional subspace of $\mathbb{C}^{2^n}$, unlike a classic [n,k,d] code which is a k dimensional subspace of $\mathbb{Z}_2^n$.*

### 2.2.3 Stabilizer Codes

Stabilizer codes are defined as a subspace that is invariant under a certain subgroup of the Pauli group. There are some nice properties of the Pauli group, resulting from Theorem 2.1, that will be useful when working with these codes.

- For every $E \in G_n$ $E^{-1} = E^\dagger$, i.e. E is unitary.

- $E^2 = I$ if the number of $\bar{Y}$ in $E$ is even, and $E^2 = -I$ otherwise.

- Since $X, Z, \bar{Y}$ all anti-commute, for every $M, N \in G_n$ $MN = -NM$ if the number of coordinates for which $N, M$ disagree and neither is the identity is odd, and $MN = NM$ if it is even.

We will represent elements of $\bar{G}_n$ as vectors in $\mathbb{Z}_2^{2n}$ by $i(E) = (\alpha, \beta)$ where $\alpha \in \mathbb{Z}_2^n$ represents the coordinates in the tensor product where E is X, $\beta$ represents the coordinates where it is Z and if both are 1, it will represent $\bar{Y} = ZX$. For example, if $E = X \otimes Y \otimes X \otimes Z \otimes I$, then $\alpha = (1,1,1,0,0)$ and $\beta = (0,1,0,1,0)$, so $i(E) = (1,1,1,0,0,0,1,0,1,0)$. It is important to notice that $i$ is a group isomorphism $i(ab) = i(a) + i(b)$.

**Definition 2.3.** *We will define $(\alpha, \beta) * (\alpha', \beta') = \alpha \cdot \beta' + \alpha' \cdot \beta (mod 2)$ as the twisted product. If we define the block matrix $\Lambda = \begin{pmatrix} 0 & I \\ I & 0 \end{pmatrix}$ then $(\alpha, \beta) * (\alpha', \beta') = (\alpha, \beta) \Lambda (\alpha', \beta')^T$*

Since $g, g' \in G_n$ commute iff they vary on an even number of coordinates, in which neither is the identity, one can see that $g, g'$ commute iff $i(g) * i(g') = 0$.

**lemma 2.1.** *Let $S = <g_1, \ldots, g_m>$ be the subgroup of $G_n$ generated by $g_1, \ldots, g_m$ and assume $-I \notin S$. The generators $g_1, \ldots, g_m$ are independent as generators (i.e. no subset of them generates $S$) iff $i(g_1), \ldots, i(g_m)$ are linearly independent.*

**Proof.** *If $i(g_1), \ldots, i(g_m)$ are linearly dependent then for some $j$, $i(g_j) = \sum_{k<j} a_k i(g_k)$ when $a_k \in \{0, 1\}$ not all 0. We can rearrange the $\{g_1, \ldots, g_m\}$ so that $i(g_l) = \sum_{k<l} i(g_k)$, and since $i$ is an isomorphism from $\bar{G}_n$ we will get that $g_l = \pm \prod_{k<l} g_k$. We will use the fact that $-I \notin S$ to conclude that $g_l = \prod_{k<l} g_k$, a contradiction to the fact that $\{g_1, \ldots, g_m\}$ is an independent set of generators.*

**lemma 2.2.** *If $i(g_1), \ldots, i(g_n)$ are independent then there is an operator $g$ such that $g$ and $g_i$ commute for $i < n$, and $g, g_n$ anti-commute.*

**Proof.** *We define $G$ as the matrix whose rows are $i(g_1), \ldots, i(g_n)$. Since $G$ if of full rank and $\Lambda$ is invertible, there is a vector $x$ satisfying $G\Lambda x^T = (0, 0, \ldots, 0, 1)$. We can pick an operator $g$ such that $i(g) = x$, and because $i(g_j) * i(g) = i(g_j) \Lambda x^T = \delta_{jn}$ we conclude that $g$ is the desired operator.*

We will now use these two lemmas and show how we can build codes from subgroups.

**Theorem 2.2.** *Let S be a subgroup of $G_n$ generated independently by $\{g_1, \ldots, g_k\}$ and $-I \notin S$. We define $\mathcal{H}_S$ as the subspace of all vectors $|\psi\rangle$ in the n qubit space fixed by $S$, i.e. $\forall g \in S : g|\psi\rangle = |\psi\rangle$. Then $S$ is an abelian group and $\dim(\mathcal{H}_S) = 2^{n-k}$.*

**Proof.** *The fact that $S$ is an abelian group is trivial, since every two elements in $G_n$ either commute or anti-commute and $-I \notin S$. Let $A \in S$ be any matrix different from the identity. For any $A \in G_n$, $A^2 = \pm I$ but since $-I \notin S$ for $A \in S$ we get $A^2 = I$. Since $A$ is unitary, it is diagonalizable, with eigenvalues $\pm 1$, but since $A \neq \pm I$ both $-1$ and $1$ must be eigenvalues of $A$. We will show that both eigenspaces must be of equal dimension $2^{n-1}$. Let $B$ be a matrix in $G_n$ that anti-commutes with $A$ by lemma 2.2, then for any eigenstate $|\psi\rangle$ with eigenvalue $\lambda$, $AB|\psi\rangle = -BA|\psi\rangle = -\lambda B|\psi\rangle$. The image of every eigenspace under the (reversible) operator is in the other eigenspace, so they have the same dimension. The space stabilized by $A$ is just the $+1$ eigenspace of $A$, so this proves the theorem for k=1. We will now use induction, we assume that the space stabilized by $<g_1, \ldots, g_l>$ is* of dimension $2^{n-l}$ and we look at the space stabilized by $< g_1, \ldots, g_{l+1} >$. From lemma 2 there is an operator $g \in G_n$ that commutes with $\{g_1, \ldots, g_l\}$ and anti-commutes with $g_{l+1}$. The mutual $+1$ eigenspace of $\{g_1, \ldots, g_l\}$ is of dimension $2^{n-l}$ (induction), and is perserved by $g$ but the $+1$ and $-1$ eigenspaces of $g_{l+1}$ are interchanged so we get that the mutual $+1$ eigenspace of $\{g_1, \ldots, g_{l+1}\}$ is half the dimension of $+1$ eigenspace of $\{g_1, \ldots, g_l\}$ one and therefore of dimension $2^{n-l-1}$. □

We define our stabilizer code $\mathcal{H}_S$ as the space stabilized by a certain group $S$. The last theorem gives us an easy way to calculate the dimension of the code. All that is left is to describe the distance of the code.

**Theorem 2.3.** *Let $S$ be a subgroup of $G_n$ generated independently by $< g_1, \ldots, g_k >$ and $-I \notin S$, and let $\mathcal{H}_S$ be the space stabilized by $S$. The distance of the stabilizer code space $\mathcal{H}_S$ is the minimal weight of an operator $E \in G_n$ such that $E$ commutes with every element of $S$, but is not in $S$.*

**Proof.** *The distance of a code is the minimal weight of $E \in G_n$ such that there is not any constant $C$ such that $\langle \bar{i} | E | \bar{j} \rangle = C \delta_{ij}$, when $| \bar{i} \rangle$ is a basis to the code space. We can divide the elements of $G_n$ into three sets:*

- *The subgroup $S$*

- *All $E \in G_n$ such that there is an element $A \in S$ that anti-commutes with them.*

- *All the elements that commute with $S$, but do not belong to it.*

*In the first case, $\langle \bar{i} | E | \bar{j} \rangle = \langle \bar{i} | \bar{j} \rangle = \delta_{ij}$ since $| \bar{j} \rangle$ is fixed by $E$. In the second case, if $N \in S$ is the operator that anti-commutes with $E$, then $\langle \bar{i} | E | \bar{j} \rangle = \langle \bar{i} | EN | \bar{j} \rangle = \langle \bar{i} | -NE | \bar{j} \rangle = -\langle \bar{i} | E | \bar{j} \rangle \Rightarrow \langle \bar{i} | E | \bar{j} \rangle = 0$. In both cases the code can be corrected. In the third case, $E$ and the generators $\{g_1, \ldots, g_n\}$ of $S$ commute and are independent, so by the previous Theorem we can divide $\mathcal{H}_S$ into 2 equal dimension subspaces, the $+1$ and $-1$ eigenspaces of $E$. This means that $E$ preserves $\mathcal{H}_S$ but acts non-trivially on it, and therefore we cannot write $\langle \bar{i} | E | \bar{j} \rangle = C \delta_{ij}$ since it is not independent on $| \bar{i} \rangle$ (C will be $\pm 1$ on the different eigenspaces).*

With these two theorems we can look at a stabilizer code as $l$ vectors $r_1, \ldots, r_l$ in $\mathbb{Z}_2^{2n}$. We demand that they are independent and that the twisted product $r_i * r_j = 0$ for $i \neq j$ (i.e. the corresponding generators commute). If d is the minimal weight of a vector $e$ such that $e*r_i = 0$ but $e \notin span\{r_1, \ldots, r_l\}$ then the stabilizer code $\mathcal{H}_S$ is a $[[n, n-l, d]]$ quantum error correcting code. We will summarize the properties of stabilizer codes in the following theorem.

**Theorem 2.4.** *Let $S$ be a subgroup of $G_n$ such that $-I \notin S$, and $S$ is generated independently by $\{g_1, \ldots, g_l\}$. The stabilizer code $\mathcal{H}_S$ has parameters $[[n, n-l, d]]$, when $d$ is the minimal weight of a vector $e \in \mathbb{Z}_2^{2n}$ such that for all $1 \leq m \leq n$ the twisted product $e * i(g_m)$ is zero, but $e \notin span\{i(g_1), ..., i(g_l)\}$.*

**CSS codes**

We will describe a specific kind of stabilizer codes due to Calderbank, Shor and Steane called CSS codes. These codes can be thought of as a quantum code built out of two classical code (see [23]), but we will focus on the description as stabilizer codes. Let $V_1$ and $V_2$ be two subspaces of $\mathbb{Z}_2^n$ of dimensions $k_1$ and $k_2$ respectively. We can define a "scalar product" on $\mathbb{Z}_2^n$ by $a \cdot b = \sum a_i b_i$, although it is not a inner product, it is a non-degenerate bilinear form, and we will use the same terminology.

**Theorem.** *If $V_1$ and $V_2$ are orthogonal subspaces of $\mathbb{Z}_2^n$ of dimensions $k_1$ and $k_2$ respectively. We can build a stabilizer code with parameters $[[n, n-k_1-k_2, d]]$ when $d$ is the minimal weight of a vector orthogonal to one of the spaces and does not belong to the other.*

**Proof.** *Let $H_i$ be a matrix whose rows are a basis for $V_i$, we will look at the stabilizer code for a subgroup $S$ whose generators are represented as rows in the matrix $\begin{pmatrix} H_1 & 0 \\ 0 & H_2 \end{pmatrix}$. Even thought there are two operators in $G_n$ that correspond with every vector in $\mathbb{Z}_2^{2n}$ we can choose consistently a representative (the one without the minus sign) and get an abelian subgroup of $G_n$. To show that $\mathcal{H}_S$ is a stabilizer code we need to show that the twisted product between two rows is zero, this is equivalent to saying that $a_i \cdot b_j = 0$ when $a_i$ is a row of $H_1$ and $b_j$ is a row of $H_2$. Since the spaces $V_1$ and $V_2$ are orthogonal the twisted product is zero, and it is indeed a stabilizer code. The dimension is $n - k_1 - k_2$ because the rows are independent in $\mathbb{Z}_2^{2n}$. The distance $d$ is the minimal weight of a vector $e = (e_1|e_2)$ such that the twisted product with every row is zero but does not belong to the space generated by the rows. This means that $e_1 \cdot b_i = 0$ and $e_2 \cdot a_i = 0$. We can assume from minimality that $e_1$ or $e_2$ is zero, and we get that the desired distance.*

### 2.2.4 Homological quantum error correcting codes

With the stabilizer formalism at hand, and the CSS example we can now define codes from a manifold. These codes are a generalization of toric codes defined by Kitaev [18]. In order to define a CSS code we need to define two orthogonal subspaces $V_1, V_2$ of $\mathbb{Z}_2^n$. The code parameters are $[[n, n - dimV_1 - dimV_2, d]]$ where d is the weight of a minimal vector orthogonal to one space and does not belong to the other.

Let $(M, S)$ be a compact connected n-manifold M with a simplicial structure S. We will denote by $C_i, Z_i, B_i$ the vector space over $\mathbb{Z}_2$ generated freely by i-simplices, i-cycles and i-boundaries of $S$ respectively (we will only be interested in homology with coefficients in $\mathbb{Z}_2$). Since $\mathbb{Z}_2$ is a field, there is a natural identification between the cochains $C^i = C_i^*$ and $C_i$. Under this identification $C_i$ is not only equipped with the boundary map $\partial : C_i \to C_{i-1}$, but also with the coboundary map $\delta : C_i \to C_{i+1}$, when we define $\delta\sigma_i$ for $\sigma_i \in C_i$ to be the (i+1)-chain satisfying

$$(\delta\sigma_i, \Delta) = (\sigma_i, \partial\Delta) \qquad (2.5)$$

for all $\Delta \in C_{i+1}$.

Before we define the code, we recall the *Poincaré dual cell decomposition*. For any simplicial structure S on a compact manifold there is a dual cellular structure $S^*$ with an isomorphism $* : C_i^S \to C_{n-i}^{S^*}$, when $*$ represent the n distinct isomorphisms induced by the correspondence between the i-simplices of S and the (n-i)-simplices of $S^*$ (and it should be clear from contest which isomorphism we are referring to). Two important properties of the $*$ isomorphism is that it preserves the the inner product, and that

$$\delta(*\sigma_i) = *\partial\sigma_i. \qquad (2.6)$$

For details see [11, p. 53].

Given a compact connected manifold M with a simplicial structure S and $1 < i < n$ such that $H_i(M, \mathbb{Z}_2) \neq 0$, i.e. with non-trivial homology group in dimension $i$, we can define $V_1 = Z_i^\perp$ and $V_2 = B_i$ both subspaces of $C_i$. The space $B_i$ is a subspace of $Z_i$ and therefore $V_1$ and $V_2$ are orthogonal and this gives us a CSS code. The parameters n,k of the code are

$n = dim(C_i) = \#\{\text{i-simplices in } S\}$, and $k = n - dim(V_1) - dim(V_2) = n - (n - dim(Z_i)) - dim(B_i) = dim(H_i(M, \mathbb{Z}_2))$.

The distance d is the minimal weight of a representative of a non zero class in $V_1^\perp/V_2$ or in $V_2^\perp/V_1$ (in a specific basis). This is equal to the minimal weight of a representative of a non zero class in $Z_i/B_i$ or in $B_i^\perp/Z_i^\perp$. The space $Z_i/B_i$ is of course just $H_i(M, \mathbb{Z}_2)$ (with the basis derived from S)). To understand $B_i^\perp/Z_i^\perp$ we will turn to the dual structure. An i-chain $\Delta$ is in $B_i^\perp$ iff for all $\sigma_{i+1} \in C_{i+1}$ we have $(\Delta, \partial\sigma_{i+1}) = 0$. Using the *-isomorphism we get that

$$\Delta \in B_i^\perp \Leftrightarrow \forall \sigma_{i+1} \in C_{i+1} : (\Delta, \partial\sigma_{i+1}) = 0 \tag{2.7}$$

$$\Leftrightarrow \forall \sigma_{i+1} \in C_{i+1} : (*\Delta, \delta(*\sigma_{i+1})) = 0 \tag{2.8}$$

$$\Leftrightarrow \forall \sigma_{i+1} \in C_{i+1} : (\partial(*\Delta), *\sigma_{i+1}) = 0 \tag{2.9}$$

$$\Leftrightarrow \forall \sigma^*_{n-i-1} \in C^*_{n-i-1} : (\partial(*\Delta), \sigma_{n-i-1}) = 0 \Leftrightarrow *\Delta \in Z^*_{n-i}. \tag{2.10}$$

Where 2.7 just the definition of $B_i^\perp$, 2.8 follows from 2.6 and the fact that $*$ preserves the inner-product, 2.9 follows from 2.5, and finally 2.10 is from the fact that $*$ is an isomorphism.

This shows that $*B_i^\perp = Z^*_{n-i}$, and with the same argument it can be shown that $*Z_i^\perp = B^*_{n-i}$. The space $B_i^\perp/Z_i^\perp$ is just $H_{n-i}(M, \mathbb{Z}_2)$ <u>with the basis of dual cells</u> (it is importent to notice that in coding theory the code <u>is not invariant</u> under change of base and we are looking at a linear space with a specific base).

**Definition 2.4.** *Let M be a manifold and S be a cellular structure on it. We define the $i^{th}$ combinatorial systole of (M,S), $csys_i(M,S)$ to be the minimal weight of an i-cycle $\sigma$ not null-homological in the basis of S i-cells $csys_i(M,S) = \min\limits_{0 \neq \sigma \in Z_i^S/B_i^S} w(\sigma)$.*

We write $Z_i^S/B_i^S$ instead of $H_i(M, \mathbb{Z}_2)$ to emphasize that we are looking at this linear space with a specific base defined by S. The combinatorial systole is a discrete analog of a $\mathbb{Z}_2 - systole$ in Riemannian geometry that will be defined later on. The next theorem summarizes the results we had for homological codes.

**Theorem 2.5.** *Let (M,S) be a n-dimensional compact manifold M, and a simplicial structure S on it. If $H_i(M, \mathbb{Z}_2) \neq 0$ then we can define a code* *with parameters* $[[ |S_i|, dim(H_i(M, \mathbb{Z}_2)), min\{csys_i(M, S), csys_{n-i}(M, S^*)\}]]$, *when $S_i$ is the set of i-simplices in S and $S^*$ is the dual structure.*

## 2.3 The Conjecture

Two parameters that measure the quality of a [[n,k,d]] code are the rate $R = \frac{k}{n}$ and the relative distance $\delta = \frac{d}{n}$. The rate measures how efficiently the code transfers information, and the relative distance measures how fault-tolerant the code is. Asymptotically good codes, or plainly good codes, is a family of codes such that both $R$ and $\delta$ are bounded away from zero as $n$ goes to infinity. It has been proven that good codes exist [5], but one can ask: does good homological codes exist? The answer to this question is, to the best of my knowledge, still open.

In his paper [26] G. zémor showed that for some examples of homological codes we get that

$$R\delta^2 \leq C \cdot n^{-2} \tag{2.11}$$

when $C$ is a constant, and conjectured that the inequality 2.11 holds for all homological codes. If 2.11 holds, it is obvious that they are no good homological codes.

A variation of 2.11, that looks only at the distance, is

$$d^2 \leq C \cdot n \tag{2.12}$$

If we look at a family of simplicial structures on a constant topological manifold, then $k = dim(H_i(M, \mathbb{Z}_2))$ is constant, and both inequalities are equivalent. This equivalence is not true in the general case. There is a strong connection between the distance of the code and the geometry property of systolic freedom, therefore the bound on the distance is unrelated to the rate and the inequality 2.12 has a more visible geometric interpretation, as will be clarified later on in this thesis. This is why the inequality $d^2 \leq C \cdot n$ will be the focal point of this work.

In their paper [21], M. Freedman, D. Meyer and F. Luo show a counterexample to 2.12 , and this counterexample also negates zémor's conjecture 2.11 (this counterexample predates zémor's paper). This counterexample however does not prove that good codes exist, since it only shows that there exist codes with $d^2 = \Omega(n \cdot log^{1/2}(n))$ (see notation section). The work in [21] shows a strong and interesting connection between systolic geometry and homological error correcting.

Let M be a Riemannian manifold. The Riemannian structure on M enables us to talk about an area of a smooth i-simplex. We can define the area of a i-chain as $Area\left(\sum n_i \sigma_i\right) = \sum |n_i| Area(\sigma_i)$.

**Definition 2.5.** *Let M be a n-dimensional Riemannian manifold, $1 \leq i \leq n$ and assume $H_i(M, \mathbb{Z}) \neq 0$. Define the $i^{th}$ systole as the minimal area of a representative of a non-trivial element of the $i^{th}$ smooth homology group, $sys_i(M) = \inf\limits_{0 \neq \sigma \in H_i(M,\mathbb{Z})} Area(\sigma)$ (the area of a homology class is not well defined, and we are looking at a specific representative), i.e. the minimal area of a i-cycle (with coefficients in $\mathbb{Z}$) not null-homological. The definition of the $\mathbb{Z}_2$-systole is the same, only with the* inf *going over the homology with coefficients in $\mathbb{Z}_2$, $sys_i(M, \mathbb{Z}_2) = \inf\limits_{0 \neq \sigma \in H_i(M,\mathbb{Z}_2)} Area(\sigma)$.*

In dimension 2 there is a systolic inequality, also known as $\mathbb{Z}_2$ systolic rigidity [13]

**Theorem 2.6.** *There exists a constant C such that for any compact surface that is not the sphere, with any Riemannian metric on it the following inequality holds*

$$sys_1(M, \mathbb{Z}_2)^2 \leq C \cdot sys_2(M, \mathbb{Z}_2) = C \cdot Area(M) \tag{2.13}$$

The negation of rigidity, i.e. the fact that such inequality does not hold is called systolic freedom.

**Definition 2.6.** *Let M be a n-dimensional smooth manifold, and $p, q \in \mathbb{N}$ satisfy $p + q = n$. We will say that M is (p,q)-free if*

$$\inf_g \frac{sys_n(M)}{sys_p(M) \cdot sys_q(M)} = 0 \tag{2.14}$$

*when the infimum runs over all Riemannian matrices g on this manifold. The dimension n will be said to have (p,q)-weak systolic freedom if*

$$\inf_{(M,g)} \frac{sys_n(M)}{sys_p(M) \cdot sys_q(M)} = 0 \tag{2.15}$$

*when the infimum runs on all compact manifolds M of dimension n, and Riemannian matrices g on them. The definition of $\mathbb{Z}_2$-systolic freedom and weak $\mathbb{Z}_2$-systolic freedom will be the same, only with $\mathbb{Z}_2$-systoles.*

In [21] the authors show an example of $\mathbb{Z}_2-$systolic freedom, both weak and strong. For the purpose of homological codes the easier case of weak freedom is sufficient. We will explain the construction in [21], and how to get codes from systolic freedom. On the other hand, we will also show that in dimension 2, because of rigidity, the inequality $d^2 \leq C \cdot n$ holds.

We can summarize the results that will be proved in the following chapters:

i. If dimension $n = 2m$ has weak $(m, m)$ $\mathbb{Z}_2$-systolic freedom, then there exists a family of homological codes from compact $n$ dimensional manifolds, such that $d^2 = \Omega(n)$

ii. Dimension 4 has $(2, 2)$ weak freedom, which gives raise to a family of codes with $d^2 = \Omega(nlog^{1/2}(n))$

iii. The systolic rigidity in dimension 2, implies that there is a constant $C$ such that the bound $d^2 \leq C \cdot n$ holds for all homological codes coming from surfaces.

# Chapter 3

# Example of Weak $\mathbb{Z}_2$-systolic Freedom

## 3.1 Gromov's Mapping Torus Example

There is a simple example of systolic freedom by Gromov ([13] sec. 4.45). This direct construction fails to produce $\mathbb{Z}_2$-systolic freedom, but the construction in [21] modifies it in order to give rise to $\mathbb{Z}_2$-systolic freedom. We will present Gromov's example, give the intuition as to why it fails in the $\mathbb{Z}_2$-systolic case, and explain the idea behind the construction in [21], before we go over it in detail.

Let *V* be a simply-connected compact n-dimensional Riemannian manifold such that $S^1$ acts freely by isometries on V. For example $S^3$ when $S^1$ acts by rotation along the Hopf fibration - we can think of $S^3$ as $\{(z,w) \in \mathbb{C}^2 : |z|^2 + |w|^2 = 1\}$. The circle $S^1$ acts freely by isometries with the action $(z,w) \mapsto (e^{i\theta}z, e^{i\theta}w)$ . We can define $\tau_k$ as the action of $\frac{2\pi}{k}$. It is obvious that $order(\tau_k) = k$, and we can scale the metric on V for each $k$, in order to get that the minimal displacement of $\tau_k$ is by $\frac{2\pi}{k}$ (we can scale differently for different values of k, as the volume of V has no impact on the freedom as we will soon see). For $\epsilon_k > 0$, define $M_k = V \times [0, \epsilon_k]/[(x,0) \sim (\tau_k x, \epsilon_k)]$ the mapping torus of $\tau_k$ (we will chose a value for $\epsilon_k$ later). All these manifolds are the same topologically so we can think of them as different Riemannian metrices on the same (n+1)-dimensional topological manifold M.

To see that M is (n,1)-free we need to estimate $sys_n(M), sys_1(M)$ and $sys_{n+1}(M)$. The only essential (n+1)-cycle is the whole manifold, ergo $sys_{n+1}(M_k) = vol(M) = vol(V) \cdot \epsilon_k$. It is also clear that since $V \times \mathbb{R}$ is the universal cover of $M_k$ we have $sys_1(M_k) \geq \min\{k \cdot \epsilon_k, \frac{2\pi}{k}\}$ - any closed loop is lifted to a path that has to go k times $\epsilon_k$ in the last coordinate alone, or be displaced by $\frac{2\pi}{k}$, at least, in $V$. If we pick $\epsilon_k = \left(\frac{2\pi}{k}\right)^2$ we get that $sys_1(M_k) \geq \frac{2\pi}{k}$. The cycle minimizing the n-dimensional area is the cross-section $V \times \{t_0\}$ for any $t_0 \in [0, \epsilon_k]$. We can prove this with the use of calibration.

**Definition 3.1.** *Let $M$ be an n-dimensional Riemannian manifold. A p-differential form $\varphi$ for $1 \leq p \leq n$ will be called a* calibration *if $\varphi$ is closed (i.e. $d\varphi = 0$) and $||\varphi|| = 1$ (p-forms operate on p-tensors, and $||\cdot||$ is the operator norm).*

**Theorem 3.1.** *Let $M$ be a compact oriented Riemannian n-manifold, and let $S$ be a p-dimensional smooth submanifold. Define $\vec{S}(x)$ to be the p-tensor defined by the tangent space of $S$ at the point $x$. If for all $x \in S$, $\varphi(\vec{S}(x)) = 1$ then $S$ has minimal area in its smooth homology class (with coefficients in $\mathbb{Z}$).*

**Proof.** *Let $T$ be a smooth cycle homological to $S$. Then there is a (p+1)-smooth chain $A$ such that $S - T = \partial A$. From Stokes Theorem $\int_{S-T} \varphi = \int_{\partial A} \varphi = \int_A d\varphi = 0$, therefore $\int_S \varphi = \int_T \varphi$. The fact that $\varphi(\vec{S}(x)) = 1$ implies that $\int_S \varphi = area(S)$, and the fact that $||\varphi|| = 1$ implies that $\int_T \varphi \leq area(T)$ and the Theorem follows.*

In Gromov's example we can use calibration with the form $\varphi = dx \wedge dy$, when $x, y$ form the orthonormal frame for $D^2$, and see that any cross-section minimizes area. This shows that $sys_n(M) = Vol(V)$. We can combine all the results and get:

$$\frac{sys_{n+1}(M_k)}{sys_1(M_k) \cdot sys_n(M_k)} \leq \frac{vol(V) \cdot \frac{(2\pi)^2}{k^2}}{vol(V) \cdot \frac{2\pi}{k}} \xrightarrow{k \to \infty} 0 \qquad (3.1)$$

This proves that M is (n,1) systolic free. If we try to use this method to get $\mathbb{Z}_2$-systolic freedom, we face a problem when we try to estimate the systole in codimension one. The method of calibration only works with $\mathbb{Z}$ coefficients, as it uses integration on manifolds and Stokes theorem, both only applicable in oriented manifold. A nice example showing that the $\mathbb{Z}_2$-systole might be smaller then the regular systole is given in [21]. This example gives some good intuition as to why this happens.

Let us present the example: we will start with the "penny" $D^2 \times [0, \epsilon]$ and identify the opposite sides with a 180° rotation, $X = D^2 \times [0, \epsilon]/(x, 0) \sim (-x, \epsilon)$ (when we think of $D^2$ as the unit disk in $\mathbb{C}$). With calibration fields, as in Gromov's construction, we can show that $D^2 \times \{t_0\}$ minimizes $sys_2(X)$ for any $t_0 \in (0, \epsilon)$, and is therefore independent of $\epsilon$.

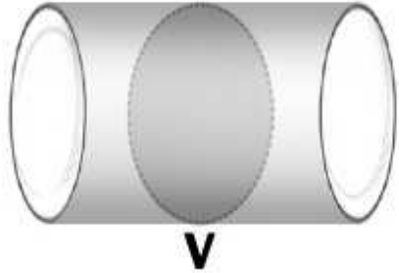

Figure 3.1: The minimal cycle with $\mathbb{Z}$ coefficients.

When we look at $\mathbb{Z}_2$-systoles, we can find smaller, not null-homological, cycles by taking out two symmetric (with regards to the rotation) "patches" and connecting them with tubes going out in opposite direction (connected by the twist), and we can make the area of this cycle as small as we want by changing $\epsilon$ and taking bigger patches.

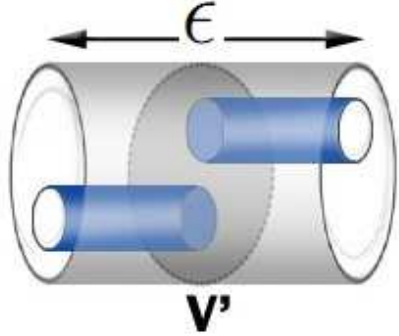

Figure 3.2: A smaller cycle with $\mathbb{Z}_2$ coefficients.

This tubular argument can be used to show that Gromov's construction with $S^3$ fails to produce $\mathbb{Z}_2$-systolic freedom. To fix this, Freedman, Meyer and Luo in [21] take instead of $S^3$ a arithmetic surface that comes from a hyperbolic space. The reason to pick an underlying hyperbolic metric, is that the isoperimetrical inequality scales linearly (as will be shown later on) and this will help us get the desired freedom.

The problem with this construction is that we will build surfaces $\Sigma_g$ with genus $g$ going to infinity, and as such the homology groups are going to have many generators. In order to keep to topology, or more specifically the homology, simple surgeries will be used to kill all the homology generators except one. We will have one free parameter (the thickness) for the surgeries, and we will show that it can be taken small enough as to have no effect on the systolic freedom.

The main reason we want only one generator for our homology is in codimension one. We have one representative we know (the cross-section, like in the previous example) and we will use the fact that every other essential cycle is homological to it, to show that it is not too small (and here is where the isoperimetric inequality will take place). When considering codes, this implies that $k = 1$ and therefore the rate goes to zero. This example does not produce good codes, but it still passes the bound given in (2.11), and is the only example of $\mathbb{Z}_2$−systolic freedom that is known (to the best of our knowledge). The question whether good homological quantum error correcting codes exist is still open, but the work in [21] is a step towards finding such codes.

## 3.2 Construction of $\Sigma_g$

The first step in building an example of $\mathbb{Z}_2-$systolic freedom will be the construction of a surface $\Sigma_g$ of genus $g$, with certain properties that will be important for systolic freedom later on. The surface $\Sigma_g$ will be a quotient of the hyperbolic plane $\mathbb{H}^2$ by a Fuchian group.

### 3.2.1 Fuchsian Groups

Let us will look at the hyperbolic plane in the upper half-plane model. In this model, the group of orientation preserving isometries $Isom^+(\mathbb{H}^2)$ is the projective special linear group, $Isom^+(\mathbb{H}^2) = PSL_2(\mathbb{R})$ acting by Möbius transformations. The trace is not well defined on $PSL_2(\mathbb{R})$, but its absolute value is. We will regard $Tr(T) = |tr(T)|$ as the trace. We want to look at subgroups $\Gamma \subset PSL_2(\mathbb{R})$ such that $\mathbb{H}^2/\Gamma$ is a manifold, a sufficient condition for this is that $\Gamma$ is discrete without fixed points.

**Definition 3.2.** *A discrete subgroup $\Gamma \subset PSL_2(\mathbb{H}^2)$ will be called a Fuchsian group.*

The fixed points of a transformation $T \in PSL_2(\mathbb{R})$ are the solutions of $z = \frac{az+b}{cz+d}$. It is easy to see that the number of distinct solutions depends on the trace of T. In the case of $Tr(T) > 2$, a straightforward calculation shows that $T$ has two distinct fixed points in $\mathbb{R} \cup \{\infty\}$.

**Definition 3.3.** *Let T be an element of $PSL_2(\mathbb{R})$. If $Tr(T) > 2$, T will be called a* hyperbolic *transformation.*

Since any two points in $\mathbb{R} \cup \{\infty\}$ are limit points of a unique geodesic in $\mathbb{H}^2$, any hyperbolic transformation T maps this unique geodesic (called the axis of T) connecting its fixed points, onto itself. It is important to notice that the fixed points are not in $\mathbb{H}^2$, but limit points of $\mathbb{H}^2$. This means that $T$ acts without fixed points on $\mathbb{H}^2$.

**Proposition 3.1.** *Let T be a hyperbolic transformation, then the translation length of T, $L(T) = \inf\limits_{x \in \mathbb{H}^2} d(x, T(x))$ satisfies $Tr(T) = 2cosh\left(\frac{L(T)}{2}\right)$.*

**Proof.** The matrix $T$ has $Tr(T) > 2$ and therefore has 2 distinct real eigenvalues, consequently it is conjugate in $SL_2(\mathbb{R})$ to a matrix of the form $\tilde{T} = \begin{pmatrix} \lambda & 0 \\ 0 & \frac{1}{\lambda} \end{pmatrix}$, $1 \neq \lambda > 0$. Conjugation by an isometry does not change the translation length or the trace, so we can look at the action of the diagonalized matrix. We will show that the minimal translation will occur on the axis (the imaginary line in this case), and prove that the equation holds. We will use two equations given in [16, Theorem 1.2.6].

$$cosh(d(z,w)) = 1 + \frac{|z-w|^2}{2Im(z)Im(w)} \tag{3.2}$$

$$cosh\left(\frac{d(z,w)}{2}\right) = \frac{|z-\bar{w}|}{2\sqrt{Im(z)Im(w)}} \tag{3.3}$$

We have $\tilde{T}(z) = \lambda^2 z$. From (3.2) we can see that if we rotate $z$ and $\tilde{T}(z) = \lambda^2 z$ to the imaginary line we maximize the denominator while keeping the numerator constant. This shows that the minimal translation occurs on the imaginary line - the axis of $\tilde{T}$. We can now use (3.3) for $z = ix$ on the axis and see that

$$cosh\left(\frac{d(z,\tilde{T}(z))}{2}\right) = \frac{|ix + i\lambda^2 x|}{2\sqrt{x\lambda^2 x}} = \frac{1}{2}\left(\lambda + \frac{1}{\lambda}\right) = \frac{Tr(T)}{2} \tag{3.4}$$

and this proves the proposition.

We can how define the Fuchsian groups we are going to work with:

**Definition 3.4.** Let $p$ be a prime number with $p \equiv 3 \pmod{4}$, and $N \geq 2$ be an integer. We define

- $\Gamma_{(p,-1)} = \left\{ A = \begin{pmatrix} a + b\sqrt{p} & -c + d\sqrt{p} \\ c + d\sqrt{p} & a - b\sqrt{p} \end{pmatrix} ; a,b,c,d \in \mathbb{Z}, det(A) = 1 \right\} / \sim$

- $\Gamma_{(p,-1)}(N) = \left\{ A \in \Gamma_{(p,-1)} : A = \begin{pmatrix} 1 + N(a+b\sqrt{p}) & N(-c+d\sqrt{p}) \\ N(c+d\sqrt{p}) & 1 + N(a-b\sqrt{p}) \end{pmatrix} \right\}$

It is important to remember that we define $\Gamma_{(p,-1)}$ and $\Gamma_{(p,-1)}(N)$ as subgroups of $PSL_2(\mathbb{R})$ not $SL_2(\mathbb{R})$.

The group $\Gamma_{(p,-1)}$ is the Fuchsian group derived from the quaternion division algebra $\frac{(p,-1)}{\mathbb{Q}}$. The condition that $p \equiv 3 \pmod 4$ guarantees that this is indeed a division algebra. The subgroup $\Gamma_{(p,-1)}(N)$ is its principle congruence subgroup with regard to N.

**Proposition 3.2.** *The subgroups $\Gamma_{(p,-1)}$ and $\Gamma_{(p,-1)}(N)$ satisfy:*

i. *The groups $\Gamma_{(p,-1)}$ and $\Gamma_{(p,-1)}(N)$ are Fuchsian groups, containing only hyperbolic elements and the unit element.*

ii. *The groups $\Gamma_{(p,-1)}$ and $\Gamma_{(p,-1)}(N)$ are cocompact subgroups of $PSL_2(\mathbb{R})$.*

iii. *The group $\Gamma_{(p,-1)}(N)$ is a normal subgroup of $\Gamma_{(p,-1)}$ with bounded index $c_1 \cdot N^2 \leq [\Gamma_{(p,-1)} : \Gamma_{(p,-1)}(N)] \leq c_2 \cdot N^3$.*

iv. *The genus $g(N) = genus(\mathbb{H}^2/\Gamma_{(p,-1)}(N))$ is proportional to the index $[\Gamma_{(p,-1)} : \Gamma_{(p,-1)}(N)]$.*

**Proof.** *Parts (i),(ii) and (iii) are general properties of Fuchsian groups derived from division quaternion algebra, and their principle congruence subgroups. The proofs of these facts can be found in [16, ch. 5]. To prove (iv) we will use the Gauss-Bonnet theorem. The area of $\mathbb{H}^2/\Gamma_{(p,-1)}(N)$ is $area(\mathbb{H}^2/\Gamma_{(p,-1)}) \cdot [\Gamma_{(p,-1)} : \Gamma_{(p,-1)}(N)]$. The manifold $\mathbb{H}^2/\Gamma_{(p,-1)}(N)$ has constant curvature $-1$ so from the Gauss-Bonnet theorem the area is proportional to the Euler characteristic, which in turn is proportional to the genus, proving the Proposition.*

We fix $p$, an we now define $\Sigma_g = \mathbb{H}^2/\Gamma_{(p,-1)}(N)$, when $g = g(N) = genus(\mathbb{H}^2/\Gamma_{(p,-1)}(N))$. This notation specifies the genus $g$ and not $N$, since the genus will be the parameter of interest when proving systolic freedom.

### 3.2.2 Geometric properties of $\Sigma_g$

We will show some geometric properties of $\Sigma_g$ that will be needed later on to prove $\mathbb{Z}_2-$systolic freedom. The first result is an immediate consequence of the Gauss-Bonnet theorem, and the fact that the only 2−cycle is the whole manifold (as holds for any connected compact manifold).

**Proposition 3.3.** *The second systole of $\Sigma_g$ satisfies $sys_2(\Sigma_g, \mathbb{Z}_2) = \Theta(g)$.*

We will now use the spectral theory of the Laplace-Beltrami operator to bound the isoperimetric constant.

**Theorem 3.2.** *There exists a constant $C_1$, such that all $\Sigma_g$ and all 2-dimensional submanifolds $A$ of $\Sigma_g$ with piecewise smooth boundary $\partial A$ satisfying $Area(A) \leq \frac{1}{2} \cdot Area(\Sigma_g)$, the isoperimetric inequality $Area(A) \leq C_1 \cdot length(\partial A)$ holds.*

**Proof.** *This theorem states that there is a uniform lower bound on the Cheeger constant of $\Sigma_g$. This is an immediate consequence of Selbeg's thorem, Cheeger's theorem and Buser's theorem found in [20, ch. 4 ], [4, ch. 8] and [3].*

We will now show a lower bound on the injectivity radius of $\Sigma_g$.

**Definition 3.5.** *Let M be a Riemannian manifold. The* injectivity radius *of M at a point $p \in M$ is the supremum of the radii for which the exponential map is injective. The injectivity radius of M is the infimum of the injectivity radius at p for all points p in M.*

The surface $\Sigma_g$ has $\mathbb{H}^2$ as its universal cover, and since $\mathbb{H}^2$ has constant negative curvature it is easy to see that the injectivity radius is half the length of the shortest closed geodesic curve (and that this curve is essential).

**Theorem 3.3.** *The injectivity radius of $\Sigma_g$ satisfies $inj.rad\,(\Sigma_g) = \Omega(log(g))$.*

**Proof.** *The hyperbolic plane $\mathbb{H}^2$ is the universal cover of $\Sigma_g$ with deck transformation group $\Gamma_{(p,-1)}(N)$, therefore from the previous remark, $inj.rad\,(\Sigma_g)$ is half the minimal translation length of $Id \neq T \in \Gamma_{(p,-1)}(N)$. We have shown a connection between the translation length and trace for hyperbolic transformations, so we will first show that for all $Id \neq T \in \Gamma_{(p,-1)}(N)$, $Tr(T) \geq N^2 - 2$.*

*Let $Id \neq T \in \Gamma_{(p,-1)}(N)$, then $Tr(T) = |2 + 2aN| \geq 2N|a| - 2$. If we write the determinant of T explicitly, we see that $det(T) = 1 + 2aN \pmod{N^2}$, however $det(T) = 1$ by definition. This proves that $N|2a$, therefore $|a| \geq \frac{N}{2}$. We can conclude that $Tr(T) \geq 2N|a| - 2 \geq N^2 - 2$.*

*With Proposition 3.1 we now have*

$$exp\left(\frac{L(T)}{2}\right) \geq cosh\left(\frac{L(T)}{2}\right) = \frac{Tr(T)}{2} \geq \frac{N^2}{2} \qquad (3.5)$$

*We can connect N to the genus using Proposition 3.2 and get*

$$exp\left(\frac{L(T)}{2}\right) \geq C \cdot g^{\frac{2}{3}} \tag{3.6}$$

*when C is some constant. If we take the log of both sides of this inequality we get that*

$$L(T) \geq \frac{4}{3}log(g) + log(C) \tag{3.7}$$

*and this proves that $inj.rad\,(\Sigma_g) = \Omega(log(g))$.* □

The next step will be to use the isoperimetric inequality in Theorem 3.2 to prove a bound on the diameter of $\Sigma_g$. For that we will need to use the smooth coarea formula. We will state here a weaker version of it, that will be sufficient for our needs. For the general form and proof see [6, ch. 3].

**The Smooth coarea formula.** *Let M be a Riemannian manifold, and $\phi : M \to \mathbb{R}$ be a smooth map, then*

$$\int_M f ||grad(\phi)|| dV = \int_{\mathbb{R}} dt \int_{\phi^{-1}(t)} f dA \tag{3.8}$$

*when $f \in L_1(M)$ and dV and dA are the Riemannian measures on M and $\phi^{-1}(t)$ respectfully.*

From Sard's theorem [22], $\phi^{-1}(t)$ is a submanifold for all but a set of measure zero, so the right hand side of the coarea formula is well defined.

**Theorem 3.4.** *Let $d_g = diameter(\Sigma_g)$. There exist a constant $C_2$, such that for all values of g, $d_g \leq C_2 \cdot log(g)$.*

**Proof.** *Let $p \in \Sigma_g$. We define $\phi_p(x) = d(p, x)$ the distance function from p, $B_p(t)$ as the geodesic ball of radius t with center at p and $\gamma_p(t) = \partial B_p(t)$ the geodesic sphere of radius t with center at p.*
*From the coarea formula (with $f \equiv 1$)*

$$\int_{B_p(t_0)} |grad(\phi_p)| dM = \int\limits_{t=0}^{t_0} length(\gamma_p(t)) dt. \tag{3.9}$$

*Since $|grad(\phi_p)| \equiv 1$, the left integral is equal to $Area(B_p)$.*
*We can differentiate both sides and get*

$$\frac{d}{dt}Area(B_p(t)) = length(\gamma_p(t)). \tag{3.10}$$

For $t$ such that $Area(B_p(t)) \leq \frac{1}{2} \cdot Area(\Sigma_g)$, we can use the isoperimetric inequality:

$$\frac{d}{dt}Area(B_p(t)) \geq \frac{1}{C_1}Area(B_p(t)). \tag{3.11}$$

At first glance, it might be possible to think that for $t$ smaller then the value needed to reach half the total area, we would get

$$Area(B_p(t)) \geq const \cdot e^{\frac{t}{C_1}}. \tag{3.12}$$

This is of course not true, because $Area(B_p(0)) = length(\gamma_p(0)) = 0$. We can, however, use the bounds from below on the area and injectivity radius of $\Sigma_g$ shown previously, to show that we can find $t_\epsilon > 0$ such that for all $\Sigma_g$ the geodesic ball with radius $t_\epsilon$ and center at p has the same area, and $\frac{1}{2}Area(\Sigma_g) > Area(B_p(t_\epsilon)) = a_\epsilon > 0$. We can say that for $t_\epsilon < t$ with $\frac{1}{2}Area(\Sigma_g) \geq Area(B_p(t))$

$$Area(B_p(t)) \geq a_\epsilon e^{\frac{t-t_\epsilon}{C_1}}. \tag{3.13}$$

We can now take $x, y \in \Sigma_g$ of maximum distance, i.e. $d_g = d(x, y)$. If we take $t_x, t_y$ to be the radii for which the area of the geodesic ball with center in $x, y$ respectfully reaches half the total area, then from the triangle inequality $t_x + t_y \geq d_g$. We can assume without lose of generality that $t_x \geq \frac{d_g}{2}$.

$$\frac{1}{2}Area(\Sigma_g) = Area(B_x(t_x)) \geq a_\epsilon e^{\frac{t_x - t_\epsilon}{C_1}} \tag{3.14}$$

The area of $\Sigma_g$ is proportional to the genus so from (3.14) it is clear that there exists a constant C such that

$$t_x \leq C \cdot log(g). \tag{3.15}$$

We can now take $C_2 = 2 \cdot C$ and we conclude that

$$d_g \leq C_2 \cdot log(g). \tag{3.16}$$

$\square$

As said previously, we want to "kill" almost all the generators of $H_1(\Sigma_2, \mathbb{Z}_2)$ with surgeries. To do this we need to find representatives with bounded length.

**Proposition 3.4.** *There exist a constant $C_3$, such that for every $\Sigma_g$ there is an embedded wedge of closed simple geodesics curves $W_1, \ldots, W_{2g}$ spanning $H_1(\Sigma_g, \mathbb{Z}_2)$, all satisfying $length(W_i) \leq C_3 \cdot log(g)$.*

**Proof.** *We can define the Dirichlet region for $\Gamma_{(p,-1)}(N)$ centered at $y \in \mathbb{H}^2$ to be the set*

$$D_y(\Gamma_{(p,-1)}(N)) = \{x \in \mathbb{H}^2 | \forall T \in \Gamma_{(p,-1)}(N), \, d(y,x) \leq d(T(y),x))\} \quad (3.17)$$

*The Dirichlet region is a fundamental domain of $\Gamma_{(p,-1)}(N)$ (it is a closed domain, with a representative of every equivalence class and this representative is unique in the interior). This means that*

$$\mathbb{H}^2/\Gamma_{(p,-1)}(N) = D_y(\Gamma_{(p,-1)}(N))/\Gamma_{(p,-1)}(N) \quad (3.18)$$

*The boundary of $D_y(\Gamma_{(p,-1)}(N))$ is a piecewise geodesic curve [16, ch. 3]. If we look at the geodesic sections of the boundary $W_1, \ldots, W_{2g}$, they are connected by vertices and it is easy to see that all the vertices are projected to the same point in $\Sigma_g$. This means that the projection of the boundary is a wedge of closed geodesic curves, and they must span $H_1(\Sigma_g, \mathbb{Z}_2)$ because $D_y(\Gamma_{(p,-1)}(N)) \setminus \partial D_y(\Gamma_{(p,-1)}(N))$ is contractible. It is also obvious that the maximum length of any such curve can not be longer then twice the diameter, finishing the proof of the proposition.*

Proposition 3.4 shows that $inj.rad(\Sigma_g) = \mathcal{O}(log(g))$, so with Theorem 3.3 we know that $inj.rad(\Sigma_g) = \Theta(log(g))$

## 3.3 Construction of the Mapping Torus $M_g$

The next step in constructing the example of $\mathbb{Z}_2$−systolic freedom, will be finding an appropriate isometry $\tau : \Sigma_g \to \Sigma_g$ and building the mapping torus $M_g = \Sigma_g \times [0,1]/[(x,0) \sim (\tau(x),1)]$. To find this isometry we will look at a sequence of $\mathbb{H}^2/\Gamma_{(p,-1)}(N_k)$, when $N_k$ goes to infinity (and consequently so does the genus). The fact that we are not looking at all the possible $\Sigma_g$, but only on a subsequence, has no effect on the weak $\mathbb{Z}_2$−systolic freedom since we only need to show that the infimum goes to zero.

**lemma 3.1.** *There exists two sequences of natural numbers $\{N_k\}_{k=1}^{\infty}$, and $\{L_k\}_{k=1}^{\infty}$ such that:*

i. *$L_k$ divides $N_k$*

ii. *If we define $g_k = genus(\mathbb{H}^2/\Gamma_{(p,-1)}(N_k))$ and $h_k = genus(\mathbb{H}^2/\Gamma_{(p,-1)}(L_k))$, then $log(g_k) = \Theta(log(h_k)^2)$*

iii. *$L_k$ and $N_k$ go to infinity (and consequently, so do $g_k$ and $h_k$)*

**Proof.** *We can define $L_k = k$, and $N_k = k^{\lfloor log(k) \rfloor}$. It is clear that $L_k | N_k$ and both series go to infinity. To finish the proof we need to prove that $log(g_k) = \Theta(log(h_k)^2)$. From Proposition 3.2 we know that*

$$log(g_k) = \Theta(log(N_k)) \tag{3.19}$$

*The same of course goes for $h_k$ and $L_k$. From this we get that*

$$log(g_k) = \Theta(log(N_k)) = \Theta(log(k^{\lfloor log(k) \rfloor})) = \Theta(log(k)^2) = \Theta(log(h_k)^2). \tag{3.20}$$

We will use $h_k$ to build a isometry that has a large order and a large translation length.

**Proposition 3.5.** *Let $g_k, N_k, h_k$ and $L_k$ be as in Lemma 3.1. There exists an isometry $\tau_k : \Sigma_{g_k} \to \Sigma_{g_k}$ such that:*

i. *$order(\tau_k) = \Omega(log^{1/2}(g_k))$*

ii. *There exists a normal Fuchsian subgroup $H_k \triangleleft Isom^+(\mathbb{H})$ such that $\tau_k \in H_k$, $\Gamma_{(p,-1)}(N_k) \triangleleft H_k$ and $inj.rad(\mathbb{H}/H_k) = \Omega(log^{1/2}(g_k))$*

**Proof.** *From Lemma 3.1, $L_k | N_k$. This means that $\Gamma_{(p,-1)}(N_k) \triangleleft \Gamma_{(p,-1)}(L_k)$, and therefore $\Sigma_{g_k} = \mathbb{H}^2/\Gamma_{(p,-1)}(N_k)$ is a covering space of $\Sigma_{h_k} = \mathbb{H}^2/\Gamma_{(p,-1)}(L_k)$. Define $\pi : \Sigma_{g_k} \to \Sigma_{h_k}$ to be the projection map. Let $\alpha$ be the shortest essential closed simple curve in $\Sigma_{h_k}$, then from Theorem 3.3, $length(\alpha) = \Omega(log(h_k))$. We can pick a base point $q \in \Sigma_{g_k}$ for $\alpha$, and lift $\alpha$ to $\tilde{\alpha}$, the unique lift with base point $q$ in $\Sigma_{g_k}$.*

*Define $\tau_k$ as the unique deck transformation (and isometry) that connects both ends of $\tilde{\alpha}$. Let $n = order(\tau_{g_k})$ then if we concatenate $\tilde{\alpha}, \tau_{g_k}\tilde{\alpha}, \ldots, \tau_{g_k}^{n-1}\tilde{\alpha}$* *we get a closed geodesic loop. This means that $n \cdot length(\tilde{\alpha}) \geq inj.rad(\Sigma_{g_k})$. From Proposition 3.4 we know that*

$$inj.rad(\Sigma_{g_k}) = \mathcal{O}(log(g_k)). \qquad (3.21)$$

*From Proposition 3.4 and Lemma 3.1, we get that*

$$length(\tilde{\alpha}) = length(\alpha) = \mathcal{O}(log(h_k)) = \mathcal{O}(log^{1/2}(g_k)) \qquad (3.22)$$

*We can conclude that $n = order(\tau_{g_k}) = \Omega(log^{1/2}(g_k))$. This fact should be intuitively clear, since if the length of the path is $\mathcal{O}(log^{1/2}(g_k))$ we will have to go at least $\Omega(log^{1/2}(g_k))$ times before we go a length of $\Omega(log(g_k))$, the needed length for a geodesic curve to be closed.*

*The proof of the second part of the theorem is immediate, since we can define $H_k = \Gamma_{(p,-1)}(L_k)$ and we have already shown that it satisfies all the demands of the theorem.* □

In order not to make the notation too cumbersome, we will neglect the $k$ index. We will work with $\Sigma_g$ and assume that the isometry $\tau : \Sigma_g \to \Sigma_g$ exists, remembering that we are not working with all the possible genera (and surfaces), but at a subsequence.

With the mapping $\tau : \Sigma_g \to \Sigma_g$ at hand we can define $M_g = \Sigma_g \times [0,1]/[(x,0) \sim (\tau(x),1)]$, the mapping torus. We have the product metric on $\Sigma_g \times [0,1]$, and since $(x,r) \to (\tau x, r+1)$ is an isomorphism in this metric, this gives a Riemannian metric on $M_g$ so that $Vol(M_g) = Vol(\Sigma_g)$.

## 3.4  Construction of $P_g$

The last part of constructing the example in [21] of $\mathbb{Z}_2$−systolic freedom, is performing surgeries on the mapping torus $M_g$ constructed earlier, to produce $P_g$ - a manifold $\mathbb{Z}_2$−homologicaly equivalent to $S^1 \times S^2$, i.e. for all $i \in \mathbb{N}$ $H_i(P_g, \mathbb{Z}_2) \simeq H_i(S^1 \times S^2, \mathbb{Z}_2)$.

Before constructing $P_g$ we will prove some easy facts about tubes in $M_g$

**Definition 3.6.** *Let $M$ be a Riemannian manifold, $\gamma$ a closed simple geodesic in $M$ and $0 < \epsilon \in \mathbb{R}$. Define the tube $T_\gamma(\epsilon)$ as $T_\gamma(\epsilon) = \{x \in M : d(x,\gamma) < \epsilon\}$.*

**Proposition 3.6.** *Let $W$ be a closed simple non-contractible geodesic curve in $\Sigma_g$. For some $z_0 \in [0,1]$ we define $\gamma = W \times z_0$, a closed simple non-contractible geodesic curve in $M_g$. For $\epsilon < \frac{1}{2}inj.rad(M_g)$ the tube $T_\gamma(\epsilon)$ is topological equivalent to a full torus. The tube is isometric to $D^2(\epsilon) \times [0, length(\gamma)]/\sim$ (identifying $t = 0$ with $t = length(\gamma)$) with the metric $ds^2 = dr^2 + r^2 d\theta^2 + cosh^2(rcos(\theta))dt^2$.*

**Proof.** *We can lift $\gamma$ to $\tilde{\gamma}$ in the universal cover $\mathbb{H}^2 \times \mathbb{R}$. The metric around $\tilde{\gamma}$ in Fermi coordinates [4, Ch. 1] in the $\mathbb{H}^2$ component, and in the standard coordinated in the $\mathbb{R}$ component is*

$$ds^2 = d\rho^2 + cosh^2(\rho)dt^2 + dz^2 \tag{3.23}$$

*Where $t$ is the arclength parameter along the geodesic $\tilde{\gamma}$ starting at some $x_0 \in \tilde{\gamma}$, $\rho$ is the distance of the $\mathbb{H}^2$ component from the geodesic, and $z$ the standard coordinate for $\mathbb{R}$. Since $\epsilon < \frac{1}{2}inj.red(M_g)$ the projection to $M_g$ would just be the identification of $t = 0$ with $t = length(\tilde{\gamma}) = length(\gamma)$. We can now switch to $(r, \theta, t)$ coordinates where $z = z_0 + r sin(\theta)$, $\rho = r cos(\theta)$ to get the desired metric*

$$ds^2 = dr^2 + r^2 d\theta^2 + cosh^2(r\cos(\theta))dt^2 \tag{3.24}$$

*We notice that in this coordinates, $r$ is the distance from the geodesic.* □

Since we are interested in performing surgeries on a tube of small $\epsilon$ (as to not change the systolic properties too drastically), we will always assume that the $\epsilon$ chosen will be smaller then half the injectivity radius.

We will now show how one can "kill" a generator of the first homology group by surgeries.

**Proposition 3.7.** *Let $M$ be an arbitrary Riemannian manifold, and let $W$ be a non zero element of $H_1(M, \mathbb{Z}_2)$, and perform surgery on the tube based on $W$, then the first homology group of the resulting space $M'$ will have smaller dimension $dim(H_1(M', \mathbb{Z}_2)) < dim(H_1(M', \mathbb{Z}_2))$, with $W$ null-homological in $M'$.*

**Proof.** *Performing a surgery on a tube $S^1 \times D^2$ means removing $S^1 \times D^2$ and gluing $D^2 \times S^1$ on their mutual boundary $S^1 \times S^1$. Let us call the full torus taken out $T$ and the one replacing it $T'$. We will use first the Mayer-Vietoris* sequance [14, ch. 2] with $X = M$, $A = \overline{M-T}$, and $B = T$ to bound the dimension of $H_1(\overline{M-T}, \mathbb{Z}_2)$.

$$\to H_2(X, \mathbb{Z}_2) \xrightarrow{\partial_*} H_1(A \cap B, \mathbb{Z}_2) \xrightarrow{(i_*, j_*)} H_1(A, \mathbb{Z}_2) \oplus H_1(B, \mathbb{Z}_2) \to H_1(X, \mathbb{Z}_2) \to 0 \tag{3.25}$$

When $i, j$ are the inclusion maps into $A, B$ respectively.

If we define $dim(H_1(A, \mathbb{Z}_2)) = k$ and $dim(H_1(X, \mathbb{Z}_2)) = n$, then from the exact sequence it is clear that $n = k + 1 - dim(Im(i_*, j_*)) = k - 1 + dim(ker(i_*, j_*))$. This is true since the first homology of a full torus is $\mathbb{Z}_2$, and the intersection $A \cap B$ is $S^1 \times S^1$ with a first homology group $\mathbb{Z}_2^2$. We claim that $k \leq n$, this is equivalent to saying that $dim(ker(i_*, j_*)) \leq 1$, or that $(i_*, j_*)$ is not the zero transformation. The fact that this transformation is not the zero transformation is obvious, since the inclusion map $j_*$ into the full torus is not zero for the meridian.

We have shown that by taking a full torus $T$ out, we do not add any new generators to the first homology. We will now show that by inserting the torus $T'$ we "kill" one generator. Let us use the Mayer-Vietoris sequance again, this time with $X = M' = \overline{M-T} \cup_{\partial T} T'$, $A = \overline{M-T}$ and $B = T'$. The sequence is again

$$\to H_2(X, \mathbb{Z}_2) \xrightarrow{\partial_*} H_1(A \cap B, \mathbb{Z}_2) \xrightarrow{(i_*, j_*)} H_1(A, \mathbb{Z}_2) \oplus H_1(B, \mathbb{Z}_2) \to H_1(X, \mathbb{Z}_2) \to 0 \tag{3.26}$$

We are interested in $dim(H_1(X, \mathbb{Z}_2)) = k + 1 - dim(Im((i_*, j_*)))$. We want to prove that $dim(H_1(X, \mathbb{Z}_2)) < n$. With the previous proof that $k \leq n$ at hand, we need to show that $dim(Im((i_*, j_*))) = 2$. This is easy to see, since the image of the meridian here is not zero in the first coordinate (since W the axis of the tube is not null-homological in M), and zero in the second coordinate, i.e. in $T'$. On the other hand the image of the longitude is non-zero in the second coordinate showing that the dimension is two.

We have shown that by performing surgery, the dimension decreases, and it is clear that the loop on with the surgery is done, is now contractible, and therefore is null-homological after the surgery. $\square$

**Proposition 3.8.** *Let $W_1, \ldots, W_{2g}$ be a basis of $H_1(\Sigma_g, \mathbb{Z}_2)$ given in Proposition 3.4. There exist a subset (after reordering) $\{W_1, \ldots, W_k\}$ such that performing surgeries on $\{W_1 \times \{\frac{1}{2} + \frac{1}{2k}\}, W_2 \times \{\frac{1}{2} + \frac{2}{2k}\} \ldots, W_k \times \{1\}\}$ in $M_g$, results in a manifold $P_g$ that has the same homology as $S^1 \times S^2$.*

**Proof.** *We will use a version of the Mayer-Vietoris Sequence for mapping toruses given in [14, p. 151], giving the long exact sequence*

$$\cdots \to H_1(\Sigma_g, \mathbb{Z}_2) \xrightarrow{id_* - \tau_*} H_1(\Sigma_g, \mathbb{Z}_2) \xrightarrow{i_*} H_1(M_g, \mathbb{Z}_2) \to \mathbb{Z}_2 \to 0 \quad (3.27)$$

*this is equivalent to*

$$0 \to coker(id_* - \tau_*) \to H_1(M_g, \mathbb{Z}_2) \to \mathbb{Z}_2 \to 0 \quad (3.28)$$

*Let $r = rank(id_* - \tau_*) < 2g$ (if $r = 2g$ then no surgery is needed), then at most $r$ of $W_1, \ldots, W_{2g}$ lie in $Im(id_* - \tau_*)$. Let $k = 2g - r$, then we shall reorder the indexes so that $\{W_1, \ldots, W_k\}$ span $coker(id_* - \tau_*)$. There is an $\epsilon$ small enough so that the tubes $T_{\gamma_i}(\epsilon)$ in $M_g$ around $\gamma_i = W_i \times \{\frac{1}{2} + \frac{i}{2k}\}$ are distinct and have the metric given by Proposition 3.6. We can perform surgeries on $T_{\gamma_1}(\epsilon), \ldots, T_{\gamma_r}(\epsilon)$, from Proposition 3.7. This kills $coker(id_* - \tau_*)$ and therefore the resulting manifold $P_g$ has $H_1(P_g, \mathbb{Z}_2) = \mathbb{Z}_2$. Since $P_g$ is a 3-manifold, it is clear that $H_i(P_g, \mathbb{Z}_2) = \mathbb{Z}_2$ for $i = 0, 3$ and $H_j(P_g, \mathbb{Z}_2) = 0$ for $j > 3$. From Poincaré duality it is also easy to see that $H_2(P_g, \mathbb{Z}_2) = \mathbb{Z}_2$ proving that $P_g$ has the same homology groups as $S^1 \times S^2$.* □

It is important to note that the space resulting from the surgeries is a piecewise smooth manifold. Before we can talk about systolic properties we need to define a smooth structure and a Riemannian metric on $P_g$. For that we will first find a suitable metric on $T'_{\gamma_i}(\epsilon)$ so that the two metrics on $\partial T_{\gamma_i}(\epsilon)$ and on $\partial T'_{\gamma_i}(\epsilon)$ are isometric. We will then use the fact that the boundary on which the surgery was performed has a collar neighborhood, to show how they can be glued in smoothly.

We will define $T'_{\gamma_i}(\epsilon)$ as $D^2(R) \times [0, 2\pi\epsilon]$ when $R = \frac{length(\gamma_i)}{2\pi}$ with the metric:

$$ds^2 = dr^2 + r^2 \cdot cosh^2(\epsilon cos(t))d\theta^2 + \epsilon^2 dt^2 \quad (3.29)$$

It is clear that if we look at the metric on $\partial T_{\gamma_i}(\epsilon)$ (i.e. putting $r = \epsilon$), it is isometric with $\partial T'_{\gamma_i}(\epsilon)$ (i.e. putting $r = R$) by rescaling the angle $\varphi = R \cdot \theta$ and swapping $\varphi$ and $t$.

The two properties of $T'_{\gamma_i}(\epsilon)$ needed are that the boundary will be isometric to that of $T_{\gamma_i}(\epsilon)$, and that the volume is bounded. We will now show the bound on the volume.

**Proposition 3.9.** *For $\epsilon < 1$ (as will always be the case) the volume of $T'_{\gamma_i}(\epsilon)$ satisfies $Vol(T'_{\gamma_i}(\epsilon)) = \mathcal{O}(\epsilon \cdot log^2(g))$.*

**Proof.** *The volume is just the integral*

$$Vol(T'_{\gamma_i}(\epsilon)) = \int_0^{2\pi} d\theta \int_0^R dr \int_0^{2\pi\epsilon} \epsilon \cdot r \cdot cosh(\epsilon cos(t))\epsilon dt \tag{3.30}$$

*Since t is bounded by $2\pi\epsilon$ and $\epsilon$ is bounded by 1, we can bound $\epsilon \cdot cosh(\epsilon cos(t))$ by some C. From here it is clear that*

$$Vol(T'_{\gamma_i}(\epsilon)) \leq C \cdot 2\pi\epsilon \cdot length^2(\gamma_i) = \mathcal{O}(\epsilon \cdot log^2(g)) \tag{3.31}$$

*When the last inequality if because $length(\gamma_i) = \mathcal{O}(log(g))$ (Theorem 3.4).*

We will now show how to glue the metrices smoothly without effecting systolic freedom.

**Definition 3.7.** *Let M be a manifold with boundary. A collar on M is an embedding $f : \partial M \times [0, 1) \to M$.*

It can be shown that every manifold with boundary has a collar. Using the collar neighborhood, it can be proven (and easy to convince yourself geometrically) that for any two smooth manifolds $M, N$ with diffeomorphic boundaries $f : \partial M \to \partial N$ there is a unique (up to diffeomorphism) smooth structure on $M \cup_f N$. The proof of this fact is found in [15, ch. 4]. In this specific case, from the way the metric is written it is clear that the embedding of the collar can be done isometrically.

This shows that $P_g$ has a smooth structure, and since the metrices agree on the boundaries, $P_g$ has a metric that is continues, and piecewise smooth. We can use the collar neighborhood to "smoothen" the metric. Let $\{g_{ij}^k\}_{k=1}^2$ be the two metrices defined on opposing the collar neighborhoods $\partial T_{\gamma_k}(\epsilon) \times [0,1)$ and $\partial T_{\gamma_k}(\epsilon) \times (-1, 0]$. We can extend each of them smoothly to $\partial T_{\gamma_k}(\epsilon) \times (-1, 1)$ and then define the metric

$$g_{ij}(s,r) = \begin{cases} g_{ij}^1(s,r) & \text{if } r \geq 0 \\ g_{ij}^2(s,r) & \text{if } r \leq 0 \end{cases} \tag{3.32}$$

Where $r$ is the coordinate along the collar. This is the continues, piecewise smooth Riemannian metric on $P_g$ that we wish to smoothen without changing systolic properties.

Let $\phi : \mathbb{R} \to [0, 1]$ be a smooth function that satisfies

i. $\forall x < -1$, $\phi(x) = 0$

ii. $\forall x > 1$, $\phi(x) = 1$

iii. $\forall x$, $\phi(x) = 1 - \phi(-x)$

Such a function is not hard to construct explicitly. We will use $\phi$ to smoothen $g_{ij}$ by

$$g_{ij}^\delta(s,r) = \phi\left(\frac{r}{l}\right) g_{ij}^1(s,r) + \phi\left(-\frac{r}{l}\right) g_{ij}^2(s,r) \tag{3.33}$$

It is clear that $g_{ij}^\delta$ converges to $g_{ij}$ as $\delta \to 0$. We will show that it converges uniformly.

**lemma 3.2.** *The smooth metrices $g_{ij}^\delta$ converge uniformly to the piecewise smooth metric $g_{ij}$.*

**Proof.** *We need to show that $\sup |g_{ij}^\delta(s,r) - g_{ij}(s,r)| \to 0$. Assume without loss of generality $r > 0$, then $|g_{ij}^\delta(s,r) - g_{ij}(s,r)| = |\phi\left(\frac{r}{l}\right) g_{ij}^1(s,r) + \phi\left(-\frac{r}{l}\right) g_{ij}^2(s,r) - g_{ij}^1(s,r)|$ from property (iii) of $\phi$ we get $g_{ij}^1 = \phi\left(\frac{r}{l}\right) g_{ij}^1 + \phi\left(-\frac{r}{l}\right) g_{ij}^1$. We have now $|g_{ij}^\delta(s,r) - g_{ij}(s,r)| = |\phi\left(-\frac{r}{l}\right)| \cdot |g_{ij}^2(s,r) - g_{ij}^1(s,r)|$. The function $|\phi\left(-\frac{r}{l}\right)|$ is non-zero only for $r < l$, and since $g_{ij}^1, g_{ij}^2$ are uniformly continues functions that agree on $r = 0$, then the supremum of the distance goes to zero as $\delta \to 0$, and the convergence is uniformly.*

The same way we define a volume form from a Riemannian metric $dV = \sqrt{|g|}dx^1 \wedge ...dx^n$, we can define a volume form with the piecewise smooth metric $g_{ij}^\delta$ (we should notice that since the "problematic" area is covered by a single coordinate map, we can define the volume form in this coordinates specifically, without worrying about coordinate change). This gives us volumes for smooth cycles, and the ability to talk about systoles. Since $g_{ij}^\delta$ converges uniformly, for any smooth cycle the area with respect to $g_{ij}^\delta$ converges to the volume with respect to $g_{ij}$. This tells us that the systole with respect to $g_{ij}^\delta$ converges to the systole with respect to $g_{ij}$ (since we are only interested in the smallest cycle, the cycles with some bound on their length converge uniformly). In conclusion we can smoothen the piecewise smooth metric on $P_g$ without changing the systolic properties.

**A note on strong systolic freedom:** We built a series of 3-manifolds $P_g$, with different topological structures, and we will show that they prove weak systolic freedom. The Lickorish-Wallace Theorem [19] tells us that by the method of *Dehn surgeries* one can produce any 3-dimensional topological manifold from any other 3-dimensional topological manifold by a finite series of surgeries. It is possible then to perform additional surgeries in order to make the manifolds $P_g$ not only have the same homology as $S^1 \times S^2$, but to be topologically equivalent. Some extra calculations are needed to decide what bound for $\epsilon$ will result in systolic freedom, but apart from that, the proof of strong systolic freedom is very similar to the proof of weak systolic freedom. This construction is also a counterexample to a conjecture by M. Gromov, suggesting that $S^n \times S^m$ does not have $\mathbb{Z}_2$−systolic freedom for any $n, m$. The calculation of what $\epsilon$ is needed to prove strong $\mathbb{Z}_2$−systolic freedom appears in the appendix of [21].

## 3.5 Verification of systolic freedom

We are going to show that dimension three has $(2,1)$ weak $\mathbb{Z}_2$−systolic freedom, i.e. that

$$\inf_{(M,g)} \frac{sys_3(M, \mathbb{Z}_2)}{sys_1(M, \mathbb{Z}_2) \cdot sys_2(M, \mathbb{Z}_2)} = 0 \tag{3.34}$$

When the infimum runs over all compact Riemannian manifolds of dimension three.

In the previous sections we first built $\Sigma_g$ as a quotient of $\mathbb{H}^2$, then we found an isometry $\tau : \Sigma_g \to \Sigma_g$ with specific properties and defined $M_g$ as the mapping torus. The last parts of the construction was performing surgeries on tubes $T_{\gamma_i}(\epsilon)$ producing $P_g$ which have the same homology groups as $S^1 \times S^2$, and then smoothening the metric. In this construction we have one free parameter $\epsilon$ which we can make infinitesimally small. We will show that with the right choice of $\epsilon$ we get

$$\lim_{g \to \infty} \frac{sys_3(P_g, \mathbb{Z}_2)}{sys_1(P_g, \mathbb{Z}_2) \cdot sys_2(P_g, \mathbb{Z}_2)} = 0 \tag{3.35}$$

proving weak systolic freedom.

The next theorem summarizes the properties of $\Sigma_g$, $M_g$ and $P_g$ proved in the previous sections and needed to prove (3.35).

**Theorem 3.5.** *Let $\Sigma_g, \tau, M_g$ and $P_g$ be as defined in the previous sections, then the following holds:*

i. $order(\tau) = \Omega\left(log^{1/2}(g)\right)$.

ii. *There exists a normal Fuchsian subgroup $H_g \triangleleft Isom^+(\mathbb{H}^2)$ such that $\tau \in H_g$, $\Gamma_{(p,-1)}(N) \triangleleft H_g$ and $inj.rad(\mathbb{H}/H_g) = \Omega(log^{1/2}(g))$*

iii. *The surgeries are around closed simple geodesic curves $\gamma_i$, all satisfying $length(\gamma_i) = \mathcal{O}(log(g))$.*

iv. *The replacement torus $T'_{\gamma_i}(\epsilon)$ has the metric that is isometric to $T_{\gamma_i}(\epsilon)$ on their boundary, and $Vol(T'_{\gamma_i}(\epsilon)) = \mathcal{O}(\epsilon \cdot log^2(g))$.*

v. *The number of surgeries performed is bounded by $\mathcal{O}(g)$.*

vi. $Vol(\Sigma_g) = Vol(M_g) = \Theta(g)$.

vii. *There exists a constant $C$ such that for all $g$ and every $A \subset \Sigma_g$ if $Area(A) \leq \frac{1}{2}Area(\Sigma_g)$, then $Area(A) \leq C \cdot length(\partial A)$*

We will now add the demand that $\epsilon$ is small enough such that

$$\epsilon = \Theta\left(\frac{1}{g^2}\right). \tag{3.36}$$

In fact, $\epsilon = \Theta\left(\frac{1}{g^\alpha}\right)$ for any $\alpha > 1$ would work, we chose $\alpha = 2$ arbitrarily for the sake of simplicity.

To prove (3.35) we will need to estimate $sys_i(P_g, \mathbb{Z}_2)$ for $i = 1, 2, 3$.

**Theorem 3.6.** *Let $P_g$ be as previously defined, then $vol(P_g) = sys_3(P_g, \mathbb{Z}_2) = \mathcal{O}(g)$*

**Proof.** *From Theorem 3.5 (iv) we know that the area of the replacement torus satisfies $Vol(T_{\gamma_i}(\epsilon)) = \mathcal{O}(\epsilon \cdot log^2(g))$. The volume of $P_g$ is smaller then* the volume of $M_g$ plus the sum of the volumes of the replacement tori. From Theorem 3.5 (iv) (v) and (vi), and from (3.36) we get

$$Vol(P_g) \leq Vol(M_g) + \mathcal{O}(g) \cdot \mathcal{O}(\epsilon \cdot log^2(g)) \leq \mathcal{O}(g) + \mathcal{O}\left(g \cdot log^2(g) \cdot \frac{1}{g^2}\right) \tag{3.37}$$

and this proves that $Vol(P_g) = \mathcal{O}(g)$. □

We will now bound $sys_1(P_g, \mathbb{Z}_2)$:

**Theorem 3.7.** *The first systole of $P_g$ satisfies $sys_1(P_g, \mathbb{Z}_2) = \Omega(log^{1/2}(g))$.*

**Proof.** Let $\gamma$ be a essential closed, piecewise smooth loop in $P_g$. We will show that $length(\gamma) = \Omega(log^{1/2}(g))$. If $\gamma$ crosses one or more of the replacement tori, we can modify $\gamma$ to $\bar{\gamma}$ that passes only through the boundary of the tori (by taking $\bar{\gamma}$ in that part the shortest curve on the boundary between the two end points on $\partial T'_{\gamma_i}(\epsilon)$). If we look at the metric of the replacement tori $ds^2 = dr^2 + r^2 \cdot cosh^2(\epsilon cos(t))d\theta^2 + \epsilon^2 dt^2$, we see that the length of the t component does not increase when we push $\gamma$ to the boundary. For each $t = t_0$ the cross section $D^2(R) \times \{t_o\}$ has the metric $ds^2 = dr^2 + r^2 \cdot cosh^2(\epsilon cos(t_0))d\theta^2$. This is a standard disk metric times some bounded factor in the angular term. It is easy from here to convince yourself that changing $\gamma$ to $\bar{\gamma}$ can increase the length by no more then a fixed factor, meaning that $\gamma = \Omega(log^{1/2}(g))$ iff $\bar{\gamma} = \Omega(log^{1/2}(g))$. Since $\bar{\gamma}$ is disjoint from the interior of all the $T'_{\gamma_i}(\epsilon)$, we can think of it as a curve in $M_g$ of the same length. It is then suffice to show that for every essential $\gamma$ in $M_g$ we have $\gamma = \Omega(log^{1/2}(g))$.

Now, assume by negation that $length(\gamma) \neq \Omega(log^{1/2}(g))$. This means that there is a subsequence $g_k$ and appropriate curves $\gamma_{g_k}$ such that $length(\gamma_{g_k}) = o(log^{1/2}(g_k))$. As before, to avoid cumbersome notation we will omit the index, and assume $length(\gamma) = o(log^{1/2}(g))$.

We can lift $\gamma$ to a curve $\tilde{\gamma}$ in the covering space $\Sigma_g \times \mathbb{R}$. If the starting point of $\tilde{\gamma}$ is $(p, t)$ then its end point is $(\tau^d p, t + d)$ for some $d \in \mathbb{Z}$. Since this is the product metric, $|d| \leq length(\tilde{\gamma}) = length(\gamma) = o(log^{1/2}(g))$. From Theorem 3.5 (i) we have $order(\tau) = \Omega(log^{1/2}(g))$, and since $|d| = o(log^{1/2}(g))$ then for large enough $g$, $|d| < order(\tau)$. Now the fact that $d$ is smaller then the order of $\tau$ tells us that $\tau^d p \neq p$ (since $\tau^d$ has no fixed points). Both point $p$ and $\tau^d p$ in $\Sigma_g$ are projected to the same point in $\mathbb{H}^2/H_g$, this *means that their distance cannot be smaller then twice the injectivity radius of $\mathbb{H}^2/H_g$ (since the geodesic balls with radius smaller then the injectivity radius must be disjoint). From Theorem 3.5 (ii) we know that $inj.rad(\mathbb{H}^2/H_g) = \Omega(log^{1/2}(g))$ and we can conclude*

$$length(\gamma) = length(\tilde{\gamma}) \geq d(p, \tau^d p) = \Omega(log^{1/2}(g)) \qquad (3.38)$$

*which is a contradiction to the claim that $length(\gamma) = o(log^{1/2}(g))$ and concludes the proof.* □

The last part would be bounding $sys_2(P_g, \mathbb{Z}_2)$, to do this we will use the following two lemmas

**lemma 3.3.** *Let X be a 2-submanifold of $P_g$, then for almost all $t \in [0, \frac{1}{2} - \epsilon]$, X intersects $\Sigma_g \times \{t\}$ transversely.*

**Proof.** *First, we should notice that the geodesics on which the surgery was performed where confined to $\Sigma_g \times [\frac{1}{2}, 1]$, so $\Sigma_g \times [0, \frac{1}{2} - \epsilon]$ was unaffected by the surgeries. This means that $\Sigma_g \times \{t\}$ is indeed a submanifold of $P_g$ for any $t \in [0, \frac{1}{2} - \epsilon]$.*

*Let $p = (x, t) \in \Sigma_g \times \{t\} \cap X$, then the intersection is transverse in p (i.e. $T_p(\Sigma_g \times \{t\}) + T_p X = T_p P_g$) iff $\frac{\partial}{\partial t} \in T_p X$. We can look at the projection map $\pi : M_g \to [0,1]/\sim$, this gives us a map $\pi : P_g \setminus \cup_i T'_{\gamma_i}(\epsilon) \to [0,1]/\sim$. We can extend this to a smooth map $\hat{\pi} : P_g \to [0,1]/\sim$, and since our area of interest is unaffected by the surgeries, the exact extension used will be unimportant. We will look at $\hat{\pi}|_X : X \to [0,1]$. From Sard's Theorem we know that the set of critical values has measure zero. Now for $t \in [0, \frac{1}{2} - \epsilon]$ the point $p = (x,t) \in X$ is a critical point of $\hat{\pi}$ iff $\frac{\partial}{\partial t} \notin T_p X$ (since $\hat{\pi} = \pi$ the projection on the last coordinate for these values of t), This means that if X does not intersects $\Sigma_g \times \{t\}$ transversely, then t is a critical value of $\hat{\pi}$. The fact that the set of critical values has measure zero tells us that for almost all $t \in [0, \frac{1}{2} - \epsilon]$, X intersects $\Sigma_g \times \{t\}$ transversely.* □

**lemma 3.4.** *Let M be a smooth 3-dimensional compact manifold, and let X and Y be two smooth 2-submanifolds that are homological in $H_2(M, \mathbb{Z}_2)$, and that intersect transversely. Under these conditions there exists a manifold with corners $B \subseteq M$ such that $\partial B = X \cup Y$*

**Proof.** *Under these conditions, from Goresky's Theorem [10], we can triangulate $M$, such that $X$ and $Y$ are subcomplexes of the triangulation. In the triangulation, $X$ and $Y$ have no common 2-cell (since they intersect transversely). Now let the 3-chain $\sum_i \sigma_3^i$ be the homology between $X$ and $Y$ (notice that all the 3-simplexes are unique since the homology is with coefficients in $\mathbb{Z}_2$) and define $B = \cup_i \sigma_3^i$. It is clear that $\partial(\sum_i \sigma_3^i) = X - Y$ implies $\partial B = X \cup Y$ (remembering that the first $\partial$ is boundary operator on chains, and the second is the topological boundary of a set) since the topological boundary is just the union of all 2-faces that are boundaries of some $\sigma_3^i$ except that who are between two distinct simplexes.*

With these lemmas we will now bound the second systole of $P_g$.

**Theorem 3.8.** *The second systole of $P_g$ satisfies $sys_2(P_g, \mathbb{Z}_2) = \Theta(g)$.*

**Proof.** *It is clear that $sys_2(P_g, \mathbb{Z}_2) = \mathcal{O}(g)$ because*

$$sys_2(P_g, \mathbb{Z}_2) \leq Vol(\Sigma_g \times \{t_0\}) = \mathcal{O}(g) \qquad (3.39)$$

*for any $t_0 \in [0, \frac{1}{2} - \epsilon]$. Now let $X_g \subset P_g$ be the 2-cycle of minimal area. A regularity theorem due to Federer [8, ch. 5] states that in dimension lesser or equal to seven, and in codimension one, there is a cycle minimizing area and that this cycle is a smooth manifold. Therefore we can assume that $X_g$ is a smooth 2-submanifold of $P_g$. We want to prove that $Area(X_g) = \Omega(g)$, so we will assume by negation that $Area(X_g) \neq \Omega(g)$, and conclude that $Area(X_g) = o(g)$ (like for curves this is true for a subsequence, but we will regard $X_g$ as this subsequence).*

*From lemma 3.3 we know that for almost all $t \in [0, \frac{1}{2} - \epsilon]$, $X_g$ intersects $\Sigma_g \times \{t\}$ transversely. If we define for such $t$, $W_t = X_g \cap (\Sigma_g \times \{t\})$ then $W_t$ is a smooth, not necessarily connected, curve (or the empty set). We can use the coarea formula (3.8) to get*

$$o(g) = Area(X_g) \geq \int_{t=0}^{1/2-\epsilon} length(W_t) dt \qquad (3.40)$$

*This means that we can pick $0 < t_0 < t_1 < t_2 < \frac{1}{2} - \epsilon$ such that $X_g$ intersects $\Sigma_g \times \{t_j\}$ transversely and*

$$length(W_{t_j}) = o(g) \qquad (3.41)$$

hold. The 2-cycles $X_g$ and $\Sigma_g$ are homological since $H_2(P_g, \mathbb{Z}_2)$ has only one nontrivial element, so from lemma 3.4 there is a manifold with corners $B$ such that $\partial B = X_g \cup (\Sigma_g \times \{t_1\})$.

The first step towards achieving the desired contradiction, is to modify $X_g$ to some $X'_g$ such that the intersection of $X'_g$ with the interior of any of the replacement tori is empty, without changing the area drastically. That way we can think of $X'_g$ as a submanifold of $M_g$ instead of $P_g$.

We note that our only demand on $\epsilon$, so far, was that $\epsilon = \Theta(\frac{1}{g^2})$. We can show, similar to the proof of lemma 3.3, that for almost all values of $\epsilon$ (smaller then the injectivity radius), $X_g$ and $\partial T'_{\gamma_i}(\epsilon)$ intersect transversely. We will assume then, that we picked an $\epsilon$ such that they intersect transversely. The boundaries $\partial T'_{\gamma_i}(\epsilon)$ are of course null-homological, so since the intersection is a homological invariant, the intersection of $X_g$ and each $\partial T'_{\gamma_i}(\epsilon)$ is a null-homological 1-submanifold of $\partial T'_{\gamma_i}(\epsilon)$ [11, Ch. 1]. We can now define $X'_g = \left(X_g \setminus \cup_i T'_{\gamma_i}(\epsilon)\right) \cup_i \delta_i$, when $\delta_i$ is the bounding surface of $X_g \cap \partial T'_{\gamma_i}(\epsilon)$ in $\partial T'_{\gamma_i}(\epsilon)$. On each component, there are two possible bounding surfaces, and we can pick $\delta_i$ such that $\delta_i \subset B$. This is possible since $\partial B = X_g \cup \Sigma_g \times \{t_1\}$ and $\partial T'_{\gamma_i}(\epsilon) \cap \Sigma_g \times \{t_0\} = \emptyset$.

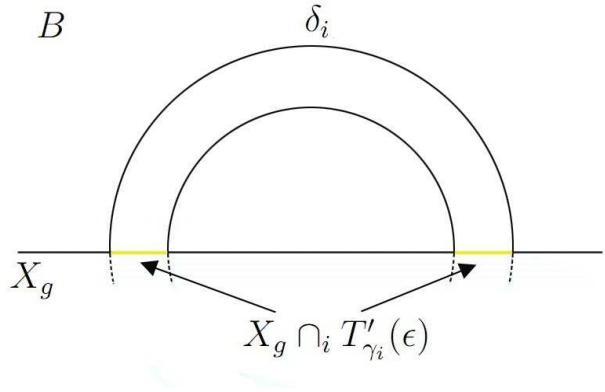

Figure 3.3: Modifying $X_g$.

The surface $X'_g$ is piecewise smooth 2-submanifold of $M_g$ (we needed smoothness to use lemmas 3.3 and 3.4, and will not need it from here on).

The area of $X'_g$ cannot be larger then the area of $X_g$ plus all the areas of $\partial T'_{\gamma_i}(\epsilon)$. For each of these tori, from (3.29)

$$Area(\partial T'_{\gamma_i}(\epsilon)) = \mathcal{O}(\epsilon log(g)) \quad (3.42)$$

so from Theorem 3.5 (v) we get

$$Area(X'_g) \leq Area(X_g) + \mathcal{O}(g) \cdot \mathcal{O}(\epsilon log(g)) < o(g) + \mathcal{O}\left(\frac{log(g)}{g}\right) = o(g) \quad (3.43)$$

We can define $B' = \overline{B \setminus \cup_i T'_{\gamma_i}(\epsilon)}$, and see that $B'$ is a $\mathbb{Z}_2$-homology between $X'_g$ and $\Sigma_g \times \{t_0\}$ in $M_g$. If we look at the boundary $\partial B'$ compared to $\partial B$ then we took out $X_g \cap T_{\gamma_i}(\epsilon)$ and added $\delta_i$ getting exactly $X'_g$, as seen in the following diagram.

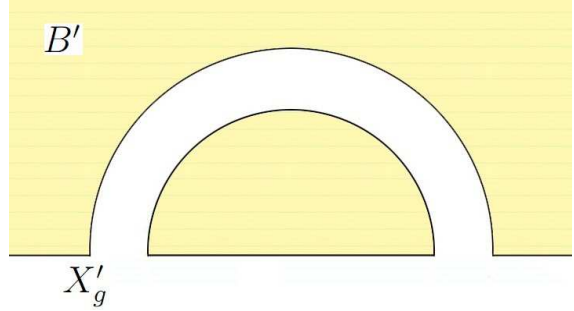

Figure 3.4: The new $B'$ and $X'_g$.

The next step will be modifying $X'_g$ to get $Z_g$ disjoint from $\Sigma_g \times \{t_1\}$. The manifolds $X'_g$ and $\Sigma_g \times \{t_j\}$ have zero intersection, then $W_{t_j}$ separates $\Sigma_g \times \{t_j\}$ into two distinct (not necessarily connected) subsurfaces. Let $Y_j$ be the one of smaller area, then from the isoperimetric inequality (Theorem 3.5 (vii)) we have

$$Area(Y_j) \leq C \cdot length(W_{t_j}) = o(g) \quad (3.44)$$

If we look at $Y_1$ and $\Sigma_g \times \{t_1\} \setminus Y_1$ then both "see" $B'$ on one side, one on the positive side ( i.e. for $t > t_1$) and one on the negative side (the sides must be different or the boundary if $B'$ would not be $\Sigma_g \times \{t_1\} \cup X'_g$.

More formally, we remember that how we constructed $B$ by triangulating $W_{t_i}$, then extending this to a triangulation of $\Sigma_g \times \{t_i\}$ and $X_g$ and then to *a triangulation of the whole manifold. We picked B to be the cycle whose boundary is $\Sigma_g \times \{t_i\} + X_g$ (thinking of both of them as cycles). Every face in $\Sigma_g \times \{t_i\}$ is either in $Y_1$ or in $Y_2$. Since this is a 3 manifold, this face is a boundary of exactly two 3-cells, one positive (i.e. with $t \geq t_i$) and the other negative. We know that this face is in $\partial B$ so one, and only one, is in B. In every connected component of $Y_1$ and $Y_2$ all the B cells must be of the same side, otherwise we could construct a curve that starts in B and ends in $B^c$ without passing through the boundary. We also know that two adjacent components of $Y_1$ and $Y_2$ must have opposite signs, otherwise assume W.L.O.G then both have B on the positive side. We can construct a curve from a point $a \in Y_1$ to a point $b \in Y_2$ that passes $X_g$ exactly ones and besides the end points is always with $t > t_i$. We can look at the faces of the cells we pass in this path, and besides the point of intersection with $X_g$, both faces must be in B (otherwise we would get a boundary that isn't in $X_g$ or $\Sigma_g \times \{t_i\}$) so we get that the two 3-cells that contain the point of $X_g$ must be in B - so it is not a boundary point of B and we get a contradiction. From here it is clear that $Y_1$ and $Y_2$ each "sees" B and on a opposite side.*

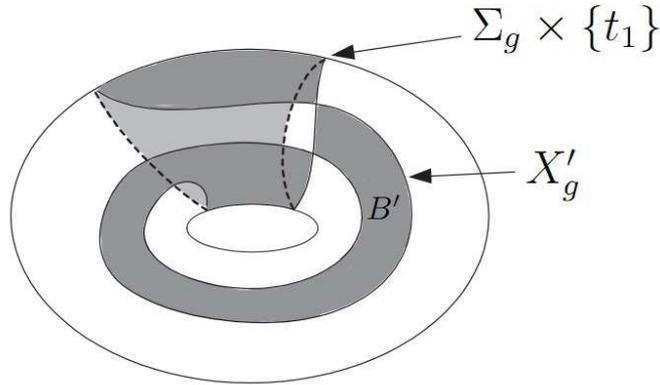

Figure 3.5: Before Modification. Image taken from [9]

*Assume without loss of generality that $Y_1$ "sees" $B'$ on the positive side. We can separate $X'_g$ from $\Sigma_g \times \{t_1\}$ by truncating $X'_g$ in $t_0, t_2$ and "glue in" $Y_2$ along $W_{t_2}$ and $\Sigma_g \times \{t_0\} \setminus Y_0$ along $W_{t_0}$ to make it a manifold without boundary (we can chose $t_0, t_1$ and $t_2$ close enough so they all "see" $B'$ on the same side). We can truncate $B'$ as well and receive $B''$. This is explained best by the following diagram.*

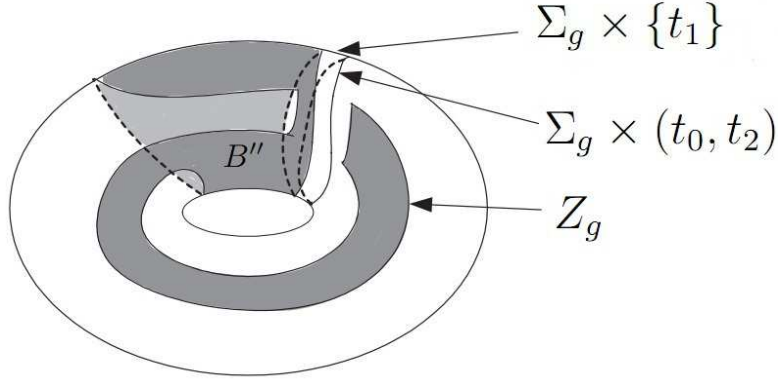

Figure 3.6: After Modification. Image taken from [9]

*Now the boarder of $B''$ is $\partial B'' = Z_g$, and*

$$B'' \subset M_g \backslash \Sigma_g \times \{t_1\} \subset M_g. \tag{3.45}$$

*The manifold $M_g \backslash \Sigma_g \times \{t_1\}$ is diffeomorphic to $\Sigma_g \times \mathbb{R}$, so it can be embedded (smooth embedding, not isometric embedding) in $\mathbb{R}^3$. From this we see that $B''$ is a 3 dimensional piecewise smooth manifold embedded in $\mathbb{R}^3$, and therefore orientable. We can use the Stoke's Theorem to arrive at the desired contradiction.*

*We can lift $Z_g$ and $B''$ to $\tilde{Z}$ and $\tilde{B}$ respectively in the covering space $\Sigma_g \times \mathbb{R}$. We will now look at the constant covector field $dt$ (where $t$ is the standard coordinate in the $\mathbb{R}$ component), and define the 2-form $\omega$ as the Hodge star dual, $\omega = *dt$ (for any k-form $\eta$, $*\eta$ is defined as the unique $(n-k)$-form such that for any k-form $\lambda$, $\lambda \wedge *\eta = <\lambda, \eta> dV_g$). If $x, y$ are coordinates for $\Sigma_g$, then it is easy to see that $\omega = *dt = \pm\sqrt{|g|}dx \wedge dy = \pm dV_\Sigma$ (the sign depends on orientation) when $dV_\Sigma$ is the volume form on the $\Sigma_g$ component. From this it is clear that $\omega$ is closed (i.e. $d\omega = 0$) and $||\omega||=1$. We can define $T = \tilde{Z} \backslash (\Sigma_g \times \{t_0\} \backslash Y_0)$ and use the Stoke's Theorem on $\tilde{B}$ to get*

$$0 = \int_{\tilde{B}} d\omega = \int_{\partial \tilde{B}} \omega = \int_{\Sigma_g \times \{t_0\} \backslash Y_0} dV_\Sigma + \int_T \omega \tag{3.46}$$

*The first integral is equal to $Area(\Sigma_g \times \{t_0\}) - Area(Y_0)$, and the second integral is smaller or equal to $Area(T)$. Since the sum is zero we have*

$$Area(\Sigma_g \times \{t_0\}) \leq Area(Y_0) + Area(T) = o(g). \tag{3.47}$$

*Where the last inequality is because $Area(Y_i) = o(g)$, and*

$$Area(T) = Area(X'_g) + Area(Y_2 = o(g) \qquad (3.48)$$

*We know that $Area(\Sigma_g \times \{t_0\}) = \Theta(g)$, arriving to the desired contradiction.* □

We can combine all the previous results are show systolic freedom

**Theorem 3.9.** *Dimension three has (2,1) weak $\mathbb{Z}_2-$systolic freedom.*

**Proof.** *We have seen that for $P_g$ we have $sys_3(P_g, \mathbb{Z}_2) = \mathcal{O}(g)$, $sys_1(P_g, \mathbb{Z}_2) = \Omega(log^{1/2}(g))$ and $sys_2(P_g, \mathbb{Z}_2) = \Theta(g)$ so*

$$\frac{sys_3(P_g, \mathbb{Z}_2)}{sys_1(P_g, \mathbb{Z}_2) \cdot sys_2(P_g, \mathbb{Z}_2)} \leq \mathcal{O}\left(\frac{g}{g \cdot log^{1/2}(g)}\right) \to 0 \qquad (3.49)$$

*proving that*

$$\inf_{(M,g)} \frac{sys_3(M, \mathbb{Z}_2)}{sys_1(M, \mathbb{Z}_2) \cdot sys_2(M\mathbb{Z}_2)} = 0 \qquad (3.50)$$

*so dimension three has (2,1) weak $\mathbb{Z}_2-$systolic freedom.* □

We can use this example of systolic freedom and get other examples of freedom by taking products.

**Theorem 3.10.** *Dimension 4 has (2,2) weak $\mathbb{Z}_2-$ systolic freedom.*

**Proof.** *This is proved in [21], and we will go over the outline of the proof here. Define $\bar{P}_g = P_g \times S^1(r)$, when $S^1(r)$ is the circle of radius $r$, with $r = g \cdot log^{-1/2}(g)$. We know that $sys_4(\bar{P}_g, \mathbb{Z}_2) = Vol(\bar{P}_g) = \Theta(g^2 log^{-1/2}(g))$.*

*We need to bound now the systole in dimension 2. From the Künneth formula [14, ch. 3]*

$$H_2(\bar{P}_g, \mathbb{Z}_2) \cong \left(H_1(P_g, \mathbb{Z}_2) \otimes H_1(S^1(r), \mathbb{Z}_2)\right) \oplus H_2(P_g, \mathbb{Z}_2). \qquad (3.51)$$

*If $[\alpha] \in H_2(\bar{P}_g, \mathbb{Z}_2) \neq 0$, then either it is non-zero in the $H_2(P_g, \mathbb{Z}_2)$ component, or it is non-zero in the $H_1(P_g, \mathbb{Z}_2) \otimes H_1(S^1(r), \mathbb{Z}_2)$ component. In the first case $Area(\alpha) \geq sys_2(P_g, \mathbb{Z}_2) = \Omega(g)$, and in the second case $Area(\alpha) \geq sys_1(P_g, \mathbb{Z}_2) \cdot sys_1(S^1(r), \mathbb{Z}_2) = \Omega(g)$. we conclude that*

$$\frac{sys_4(\bar{P}_g, \mathbb{Z}_2)}{\left(sys_2(\bar{P}_g, \mathbb{Z}_2)\right)^2} \leq \mathcal{O}\left(\frac{g^2 log^{-1/2}(g)}{g^2}\right) \to 0 \qquad (3.52)$$

We will give a precise definition to the idea of "how free" our family of manifolds is.

**Definition 3.8.** *We $S$ be a family of Riemannian n-dimensional manifolds, and $\sup_{P \in S}\{Vol(P)\} = \infty$. We say that $S$ has at least $f$ (p,q) $\mathbb{Z}_2-$systolic freedom, for some monotonically increasing function $f$, if for every $P \in S$*

$$sys_p(P, \mathbb{Z}_2) \cdot sys_q(P, \mathbb{Z}_2) = \Omega\left(f(sys_n(P, \mathbb{Z}_2)) \cdot sys_n(P, \mathbb{Z}_2)\right). \qquad (3.53)$$

We have shown that $\{P_g\}$ and $\{\bar{P}_g\}$ has $log^{1/2}$ $\mathbb{Z}_2-$systolic freedom. As we will see in the next section, stronger freedom will result in better codes, so the question whether there exists families with stronger then $log^{1/2}$ freedom is an important, yet open, question.

# Chapter 4

# Holomogical Codes and Systolic Freedom

In Theorem 2.5 we showed that if $(M, S, S^*)$ is a compact manifold $M$ with a simplicial structure $S$, a dual cellular structure $S^*$ and with $H_i(M, \mathbb{Z}_2) \neq 0$, then we can construct a quantum homological error correcting code with parameters $[[n, k, d]]$, where $n = |S_i|$ the number of $i - simplices$, $k = dim(H_i(M, \mathbb{Z}_2))$ and $d = min\{csys_i(M, S), csys_{n-i}(M, S^*)\}$. As we defined in chapter 2, $csys_i(M, S)$ is the minimal support of non-zero element of $H_i(M, \mathbb{Z}_2)$ in the basis of $S$ cells. In chapter 3 we gave in detail the construction from [21], proving that weak (and strong) $\mathbb{Z}_2$-systolic freedom exists. We will now show the connection between systolic freedom, homological codes and the $d^2 \leq C \cdot n$ bound.

## 4.1 Homological Codes from Systolic Freedom

In order to get codes from systolic freedom, we need to show that given a Riemannian manifold we can find a "nice" triangulation in respect to the metric. This is done [21, Theorem 12.8].

**Theorem 4.1.** *For any dimension $n$ there exists constants $0 < c_1, c_2, c_3 < \infty$ such that any compact n-dimensional Riemannian manifold $M$ with bounded curvature and covariant derivation of the curvature $||\nabla_h R^l_{ijk}||, ||R^l_{ijk}|| \leq 1$, and with injectivity radius at least one, there exists piecewise smooth pair*

$(S, S^*)$ of a triangulation and a dual cellulation of $M$ such that:

i. The number of cells of all dimensions in both $S$ and $S^*$ satisfies: $|S| + |S^*| < c_1 \cdot vol(M)$.

ii. For every $i$-cell $\sigma$ of $S$ or $S^*$: $vol_i(\sigma) < c_2$.

iii. For every cell of $S$ or $S^*$ the number of cells in its boundary and coboundary are less then $c_3$.

In the construction presented in the previous chapter the curvature was not bounded by one, in fact after smoothening the metric it was considerably larger then it. However, this is not a problem since we can always rescale the metric. If we rescale the metric by a factor $\alpha$, then $sys_k(M, \mathbb{Z}_2)$ rescales by a factor of $\alpha^k$ and therefore the ratio

$$\frac{sys_n(M, \mathbb{Z}_2)}{sys_q(M, \mathbb{Z}_2) \cdot sys_p(M, \mathbb{Z}_2)} \tag{4.1}$$

is scaling invariant, since $p + q = n$. As a result, we can rescale the metric so that the curvature is bounded by one (since the manifold is compact) without changing the systolic properties. We will now show how the previous theorem can be used to construct codes that pass the (2.12) bound.

**Theorem 4.2.** *Assume dimension $2m$ have (m,m) weak systolic freedom with freedom function $f$ (defined in the end of chapter 3), then there exist a family of homological quantum codes in dimension $2m$ with distance $d$ satisfying $d^2 \geq f\left(\frac{n}{c_1}\right) \cdot \frac{n}{c_1 c_2^2}$.*

**Proof.** *Let $Q^l$ be a family of compact Riemannian manifolds such that*

$$sys_m(Q^l, \mathbb{Z}_2)^2 \geq f(vol(Q^l)) \cdot vol(Q^l). \tag{4.2}$$

*The systole $sys_n(M, \mathbb{Z}_2)$ is always equal the volume in a connected manifold, which are the only manifolds we are interested in. From Theorem 4.1 there exists a triangulation $S$ and dual calculation $S^*$ such that $|S| + |S^*| < c_1 \cdot vol(Q^l)$ and that for every cell $\sigma$ we have $vol_i(\sigma) < c_2$. We can build a quantum homological code from $(Q^l, S, S^*)$ with parameters $[[n, k, d]]$ given by Theorem 2.5. We can bound $n$ by*

$$n = |S_m| < |S| + |S^*| < c_1 \cdot vol(Q^l). \tag{4.3}$$

*For every non-bounding m-cycle $\sum_j \sigma$, the volume satisfies $vol_m(\sum_j \sigma) \geq sys_m(M, \mathbb{Z}_2)$ (by the definition). This leads us to*

$$sys_m(Q^l, \mathbb{Z}_2) \leq max_{\sigma \in S_m}(vol_m(\sigma)) \cdot csys_m(M, S) \leq c_2 \cdot csys_m(M, S). \quad (4.4)$$

*When $csys_m(M, S)$ is the minimal support of a non-bounding chain in the basis defined by the elements of $S$. This inequality also holds for $csys_m(M, S^*)$ and therefore we can conclude that*

$$d^2 = min\{csys_m(M, S), csys_m(M, S^*)\}^2 \geq \left(\frac{1}{c_2}\right)^2 sys_m(Q^l, \mathbb{Z}_2)^2 \quad (4.5)$$

$$\geq f(vol(Q^l))\frac{vol(Q^l)}{c_2^2} \geq f\left(\frac{n}{c_1}\right) \cdot \frac{n}{c_1 c_2^2}. \quad (4.6)$$

$\square$

This inequality holds for for $n = 2m$ and $p = q = m$, but not for the general case $p + q = n$. For example in the case of $P_g$ with $(2,1)$ freedom, we have $d = min\{\mathcal{O}(g), \mathcal{O}(log^{1/2}(g))\} = \mathcal{O}(log^{1/2}(g)) \not\geq \mathcal{O}(g)$.

It is also important to note that part (*iii*) of Theorem 4.1 implies that this code is a LDPC (low density parity check) code and as so has a good decoding algorithm [7].

We can conclude then from the last theorem and the fact that dimension four has $(2,2)$ weak $\mathbb{Z}_2$−systolic freedom with freedom function $log^{1/2}$ (theorem 3.10) the following result is deduced in [21].

**Theorem 4.3.** *There exists a family $(\bar{P}_g, S, S^*)$ of compact four manifolds, triangulation and dual cellulation, such that the corresponding homological quantum codes are LDPC with parameters $[[n, 2, d]]$, and $d \geq \mathcal{O}(\sqrt{log^{1/2}(n) \cdot n})$.*

This is a counterexample to both the $R\delta^2 \leq \mathcal{O}(n^{-2})$ bound, and the $d^2 \leq \mathcal{O}(n)$ bound described in section 2.3.

Theorem 4.2 shows that in order to get a family of homological codes with $\delta > c > 0$ (i.e. with all the relative distances bounded away from zero), it is sufficient to prove that dimension $2m$ has $(m, m)$ freedom with freedom function that is $f(x) = \Omega(x)$. It is not yet known if this is possible to do in any dimension, and we will now show that in some cases this is the best freedom one can hope to get.

**Theorem 4.4.** *For any dimension $n = p + q$ there exists a constant $C$, such that for any compact manifold $M$ with bounded curvature and covariant derivation of the curvature $||\nabla_h R^l_{ijk}||, ||R^l_{ijk}|| \leq 1$, and with injectivity radius at least one, the following inequality holds: $sys_p(M, \mathbb{Z}_2) \cdot sys_q(M, \mathbb{Z}_2) \leq C \cdot vol^2(M)$.*

**Proof.** *This is a consequence of Theorem 4.1:*

$$sys_p(M, \mathbb{Z}_2) < |S_p| \cdot max_{\sigma \in S_P}\{vol_p(\sigma)\} < c_1 \cdot |S_p| < c_1 c_2 \cdot vol(M). \quad (4.7)$$

*This inequality also hold for $sys_q(M, \mathbb{Z}_2)$, and proves that the inequality*

$$sys_p(M, \mathbb{Z}_2) \cdot sys_p(M, \mathbb{Z}_2) \leq C \cdot vol^2(M) \quad (4.8)$$

*holds with constant $C = c_1^2 c_2^2$ for $M$. Unlike the systolic inequality, this inequality is not scaling invariant, so we cannot disregard the bound on the curvature.*

## 4.2 Bounds on Codes from Systolic Rigidity

Systolic rigidity is the negation of systolic freedom. More specifically, dimension $n$ is said to have strong $(p, q)$ $\mathbb{Z}_2$−systolic rigidity if there is a constant $C$ such that for any $n-$dimensional compact manifold $M$ and for any Riemannian metric on $M$, the following holds:

$$sys_p(M, \mathbb{Z}_2) \cdot sys_q(M, \mathbb{Z}_2) \leq C \cdot sys_n(M, \mathbb{Z}_2). \quad (4.9)$$

Strong systolic rigidity is the negation of weak systolic freedom. It has been proved [12] that dimension two has $(1,1)$ strong $\mathbb{Z}_2$−systolic rigidity. We will use the systolic inequality (4.9) to prove that there exists a constant $C$, such that for any homological quantum error correcting code that comes from a compact surface the $d^2 \leq C \cdot n$ bound holds. To do so, we will prove a Theorem dual to Theorem 4.1. Theorem 4.1 showed how to construct a "nice" triangulation of a given Riemannian manifold, we will show how given a triangulated surface we can find a "nice" Riemannian metric on it.

**Theorem 4.5.** *There exists two constants $C_1, C_2$, such that for any compact surface $M$ with a triangulation $S$ there exists a Riemannian metric $g$ on $M$ such that*

i. For any face $\sigma$ in the triangulation, $Area(\sigma) \leq C_1$

ii. $csys_1(M, S) \leq C_2 \cdot sys_1(M, \mathbb{Z}_2)$

**Proof.** *If we would give each of the triangles a flat metric of an equilateral triangle with sides of length one, we could "glue" the metrices along the sides smoothly. Unfortunately we would get singularities at the vertices (unless there are exactly six triangles meeting in that vertex). To fix this we will change the metric around the vertices.*

*Let $v$ be a vertex, and let $n$ be the number of triangles that meet at $v$. We can take the flat Euclidean metric defined before on the space, except on the disk of radius $r = \frac{1}{6}$ around each vertex. We can assume that these neighborhoods are mapped by the disk of radius $r = \frac{1}{6}$ around the origin in $\mathbb{R}^2$, and that the sector with $\frac{2\pi}{n}(i-1) \leq \theta \leq \frac{2\pi}{n}i$ belongs to the $i^{th}$ triangle.*

*We have shown previously, that if we define a continues and piecewise smooth metric we can smoothen it without changing the volume or the length of curves by more then a factor that we can make as close to one as we desire. We will show how we can define a continues piecewise smooth metric with the desired properties.*

*We will look at a vertex $v$, and as before we define $n$ as the number of triangles that meet in $v$. We will define $\eta = \frac{n}{6}$, a constant representing how the number of triangles deviates from the desired number. Since this is a triangulation, $n \geq 3$ and therefore $\eta \geq \frac{1}{2}$. We will define new coordinates $(r, \varphi)$ when $r$ is the standard radius and $\varphi = \eta \cdot \theta$. On the set $0 \leq r \leq \frac{1}{3n} = \frac{1}{18\eta}$ we will define the metric to be a standard euclidean metric*

$$ds^2 = dr^2 + r^2 d\theta^2 = dr^2 + \left(\frac{r}{\eta}\right)^2 d\varphi^2 \tag{4.10}$$

when $0 \leq \varphi \leq \eta \cdot 2\pi = n \cdot \frac{2\pi}{6}$.

On the annulus $\frac{1}{12\eta} = \frac{1}{2n} \leq r \leq \frac{1}{6}$ we will define the metric as

$$ds^2 = dr^2 + r^2 d\varphi. \tag{4.11}$$

*This metric coincides with the equilateral metric on each triangle, since each triangle has an opening of $\frac{2\pi}{6}$ radians. We are left with defining a metric*

on the annulus $\frac{1}{18\eta} = \frac{1}{3n} \leq r \leq \frac{1}{2n} = \frac{1}{12\eta}$ that coincides with the metric we previously defined when $r = \frac{1}{18\eta}, \frac{1}{12\eta}$.

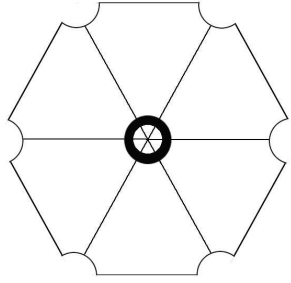

Figure 4.1: The annulus on which we "fix" the metric.

To do this we will define the metric on the annulus $\frac{1}{18\eta} \leq r \leq \frac{1}{12\eta}$ as

$$ds^2 = dr^2 + \left[\left(\frac{9\eta^2 - 4}{36\eta^3}\right) r + \left(\frac{2 - 3\eta^2}{216\eta^4}\right)\right] d\varphi^2. \tag{4.12}$$

It is easy to see that this metric is piecewise smooth and continues. We will prove that it satisfies the conditions of the theorem.

The area of each triangle is bounded by the area of an equilateral triangle with edges of length one, plus the area around the corners up to $r = \frac{1}{12\eta}$. Let us call $\sigma_1^v$ the sector with $0 \leq r \leq \frac{1}{18\eta}$, and $\sigma_2^v$ the sector with $\frac{1}{18\eta} \leq r \leq \frac{1}{12\eta}$. It is clear that

$$area(\sigma_1^v) = \int_{\frac{2\pi}{6}(i-1)}^{\frac{2\pi}{6}i} d\varphi \int_0^{\frac{1}{18\eta}} \left(\frac{r}{\eta}\right) dr = \frac{2\pi}{12 \cdot 18^2 \eta^3} < \frac{1}{77} \tag{4.13}$$

when in the last inequality we used the fact that $\eta \geq \frac{1}{2}$.

In $\sigma_2^v$ the volume form is $dV = \sqrt{\left[\left(\frac{9\eta^2-4}{36\eta^3}\right) r + \left(\frac{2-3\eta^2}{216\eta^4}\right)\right]} dr d\varphi$. The linear argument of the metric gets its maximum on the boundary, and it is equal to $\frac{1}{18^2 \eta^4}$ and $\frac{1}{12^2 \eta^2}$ at the two edges. For $\eta \geq \frac{1}{2}$ and $\frac{1}{18\eta} \leq r \leq \frac{1}{12\eta}$ it is easy to

see that
$$\left(\frac{9\eta^2 - 4}{36\eta^3}\right) r + \left(\frac{2 - 3\eta^2}{216\eta^4}\right) \leq \frac{4}{81}. \tag{4.14}$$

in consequence we can bound the volume form in $\sigma_2^v$

$$dV \leq \frac{2}{9} dr d\varphi. \tag{4.15}$$

With the bound on the volume form we can bound the volume by

$$area(\sigma_2^v) = \int_{\sigma_2^v} dV \leq \int_{\frac{1}{18\eta}}^{\frac{1}{12\eta}} dr \int_{\frac{2\pi}{6}(i-1)}^{\frac{2\pi}{6}i} \frac{2}{9} d\varphi = \frac{\pi}{486} < \frac{1}{150}. \tag{4.16}$$

We can now bound the total area of each triangle $S$ by

$$area(S) \leq \frac{\sqrt{3}}{4} + 3 \cdot \left(\frac{1}{150} + \frac{1}{77}\right) < \frac{1}{2} \tag{4.17}$$

and we have proved that all triangle areas are bounded by $C_1 = \frac{1}{2}$ finishing the first half of the proof.

We now need to show that there exists a constant $C_2$ such that

$$csys_1(M, S) \leq C_2 \cdot sys_1(M, \mathbb{Z}_2). \tag{4.18}$$

The surface $M$ is compact, and therefore there is a closed geodesic curve $\gamma$ such that $length(\gamma) = sys_1(M, \mathbb{Z}_2)$ (i.e. the infimum is achieved) [17, ch. 5]. We will show that we can find a path $G$ on the graph that is homological to $\gamma$ (with $\mathbb{Z}_2$ coefficients) and that $length(G) \leq C_2 \cdot length(\gamma)$. It is important to notice that the length of the triangles' edges in this metric are all one, so thinking about $length(G)$ in the Riemannin metric or as the distance along the underlying graph of the triangulation $S$ is the same.

we will construct the path $G$ in the following way: Let $S_0$ be the set of vertices (i.e. 0-cells of the triangulation), and we can pick $x_0 \in S_0$ such that $d(x_0, \gamma) = \min_{x \in S_0}\{d(x, \gamma)\}$. We can reparameterize $\gamma$ so that after reparameterization we have $\gamma : [0, 1] \to M$ and $d(\gamma(0), x_0) = d(\gamma, x_0)$ - the curve starts at one of the points closest to $x_0$. We can define $t_1$ as the first time when the curve $\gamma$ leaves the triangles that have $x_0$ as a vertex (it has to leave or the curve will be contractible). We can define $x_1$ as one of the vertices closest to $\gamma(t_1)$ on that edge. Since $x_1$ is one of the vertices in an edge opposing $x_0$, it is connected to $x_0$. We can continue in the same manner and get two finite sequences $x_0, ..., x_n$ and $t_0, ..., t_n$ such that

i. $0 = t_0 < t_1 < ... < t_n = 1$

ii. For all $i$, $x_i \in S_0$ and $x_n = x_0$ (we can chose$x_n$ to be $x_0$ if there is not a single option).

iii. For all $0 \leq i < n$ there is an edge $e_i$ in $S_1$ connecting $x_i$ and $x_{i+1}$.

iv. the segment $\gamma([t_{i-1}, t_{i+1}])$ is contained in the traingles that have $x_i$ as their vertex (we identify $t_0$ with $t_n$ so $t_{-1} \equiv t_{n-1}$).

We will now show that the cycle $\sum e_i$ is homological to $\gamma$ and we will show a connection between their lengths.

It is obvious from (iv) that the two paths $G$ and $\gamma$ are homotopically equivalent - we can define the homotopy from $\gamma([t_i, t_{i+1}])$ to $e_{i-1}e_{i+1}$ since there are both contained in a convex set. Homology is invariant under homotopy so $G$ and $\gamma$ are homologically equivalent.

To conclude our proof we will show that $length([t_i, t_{i+1}) \geq \frac{1}{6}$ and therefore

$$csys_1(M, S) \leq length(G) = length(\sum e_i) = n \leq 6 \cdot length(\gamma) = 6 \cdot sys_1(M, \mathbb{Z}_2). \tag{4.19}$$

In order to prove that $length([t_i, t_{i+1}]) \geq \frac{1}{6}$ we will first notice that the length of a curve going out of a vertex with a fixed angle is the same as the euclidean length (since we only changed the angular part of the metric). From this we may conclude that if the point $x$ is in a triangle that has $v$ as a vertex, and it has the coordinates $x = (r, \varphi)$ (with regard to $v$) then $d(x, v) \leq r$ (since there is a path connecting them of length $r$). If $r > \frac{1}{6}$ then we have

$$r - (\frac{1}{12\eta} - \frac{1}{18\eta}) \leq r - \frac{1}{18} < d(x, v) \leq r \tag{4.20}$$

(we just take the length on the non-euclidean part to be zero). From inequality (4.20) and the law of cosines it is easy to show that $\gamma(t_0)$ has $r < \frac{2}{3}$ (with regard to $x_0$). For all other $\gamma(t_i)$ we have $r \leq \frac{1}{2}$ with regard to $x_i$, since $x_i$ is the vertex closest to $\gamma(t_i)$ on an edge. Since $r < \frac{2}{3}$ it has to go at least a distance of $\frac{1}{6}$ to leave the triangle (considering only the euclidean part) so $length(\gamma[t_i, t_{i+1}]) \geq \frac{1}{6}$.

$\square$

With the lemma giving us a connection between a discrete triangulation and Riemannian geometry, we can finally use systolic rigidity to prove the desired square root bound on quantum homological error correcting codes, obtained from compact surfaces.

**Theorem 4.6.** *There exists a constant $C$, such that for any compact surface $M$ with triangulation $S$, the homological quantum error correcting code defined by $(M, S, S^*)$ with distance $d$ satisfies $d^2 \leq C \cdot n$.*

**Proof.** *From the previous theorem we can define a Riemannian metric such for each 2-cell $\sigma$, $area(\sigma) \leq \frac{1}{2}$ and that $csys_1(M, S) \leq 6 \cdot sys_1(M, \mathbb{Z}_2)$. From the systolic inequality we know that $sys_1(M, \mathbb{Z}_2)^2 \leq \frac{4}{3} area(M)$. Combining the results we get that*

$$csys_1(M, S)^2 \leq 36 \cdot sys_1(M, \mathbb{Z}_2) \leq 48 \cdot area(M) \leq 24 \cdot |S_2|. \quad (4.21)$$

*In a triangulation $|S_2| = \frac{3}{2}|S_1|$ (by counting boarders) so we have*

$$csys_1(M, S)^2 \leq 36 \cdot |S_1|. \quad (4.22)$$

*The number of 1-cells is the parameter $n$ of the code, and $d = min\{csys_1(M, S), csys_1(M, S^*)\} \leq csys_1(M, S)$ so we have*

$$d^2 \leq 36 \cdot n. \quad (4.23)$$

$\square$

There are examples in [26] that show codes with square root bound. This shows that the bound shown here (up to changing the constant factor) is tight.

One can try and expand this inequality to codes that come from a cellulation of a surface, not just a triangulation. This generalization is not as immediate as one might think, since the inequality $csys_1(M, S)^2 = \mathcal{O}(n)$ does not have to hold (for a counterexample just any triangulation, then divide the edges by adding more vertices and repeat), but does hold for the dual structure $S^*$. This means that one has to look at both $S$ and $S^*$. However, the work presented here makes it clear, in our opinion, that the point of interest in studying homological codes, and looking for codes ones, is in dimensions higher then two.

# Bibliography

[1] B. Podolsky A. Einstein and N. Rosen. Can quantum-mechanical description of physical reality be considered complete? *Physical Reviews*, 47, 1935.

[2] J.S. Bell. On the einstein podolsky rosen paradox. *Physics Publishing Co.*, 1:195–200, 1964.

[3] P. Buser. A note on the isoperimetric constant. *Annales Scientifiques de l'École Normale Supérieur*, 15(4):213–230, 1982.

[4] P. Buser. *Geometry and Spectra of Compact Riemann Surfaces. Birkhäuser*, 1992.

[5] A. R. Calderbank and P. W. Shor. Good quantum error-correcting codes exist. *Physical Review A*, 54:1098–1105, 1995.

[6] I. Chavel. *Riemmanian Geometry: A Modern Introduction*. Cambridge University Press, 2nd edition, 2006.

[7] G. Mitchison D.J.C. MacKay and P.L. McFadden. Sparse-graph codes for quantum error-correction. *IEEE Transactions on Information Theory*, 50:2315 – 2330, 2004.

[8] H. Federer. *Geometric Measure Theory*. Springer, 1969.

[9] M.H. Freedman. $\mathbb{Z}_2-$systolic freedom. *Geometry and Topology Monograph*, 2:113Ű12, 1999. Proceedings of the Kirbyfest.

[10] R.M. Goresky. Triangulation of stratified objects. *Proceedings of the American Mathematical Society*, 72(1), 1978.

[11] P. Griffiths and J. Harris. *Principles of Algebraic Geometry*. Wiley, 1978.

[12] M. Gromov. Filling riemannian manifolds. *Journal of Differential Geometry*, 18:1–147, 1983.

[13] M. Gromov. *Metric Structures for Riemannian and Non-Riemannian Spaces*. Springer, 1999.

[14] A. Hatcher. *Algebraic Topology*. Cambridge University Press, 2002.

[15] M.W. Hirsch. *Differential Topology*. Springer, 1994.

[16] S. Katok. *Fuchsian Groups*. The University of Chicago Press, 1992.

[17] M.G. Katz. *Systolic Geometry and Topology*. American Mathematical Society, 2000.

[18] A.Y. Kitaev. Fault-tolerant quantum computation by anyons. *Annals of Physics*, 303:2–30, 2003.

[19] R. Lickorish. A representation of oriented combinatorial 3-manifolds. *Annals of Mathematics*, 72:531–540, 1962.

[20] A. Lubotzky. *Discrete Groups, Expanding Graphs and Invariant Measures. Birkhäuser*, 1994.

[21] D.A. Meyer M.H. Freedman and F. Luo. $\mathbb{Z}_2$-systolic freedom and quantum codes. *in mathematics of quantum computation*, pages 287–320, 2002.

[22] J.W. Milnor. *Topology from the Differentiable Viewpoint*. Princeton University Press, 1965.

[23] M.A. Nielsen and I.L. Chuang. *Quantum Computation and Quantum Information*. Cambridge University Press, 2000.

[24] J. Presskil. Quantum computation. www.theory.caltech.edu/people/preskill/ph229/.

[25] P.W. Shor. Polynomial-time algorithms for prime factorization and discrete logarithms on a quantum computer. *SIAM Journal on Scientific Computing*, 26:1484–1509, 1997.

[26] G. Zémor. On cayley graphs, surface codes, and the limits of homological coding for quantum error correction. *Coding and Cryptology, second international workshop IWCC*, LNCS 5557:259–273, 2009.